\documentclass{amsart}
\usepackage[utf8]{inputenc}
\usepackage[top=1cm, bottom=2cm, left=2cm, right=2cm]{geometry}
\usepackage{amssymb}
\usepackage[cm]{fullpage}
\usepackage{wrapfig}
\usepackage{graphicx}
\usepackage{courier}
\usepackage{caption}
\usepackage{subcaption}
\usepackage{float}
\usepackage{yfonts}
\usepackage{mathrsfs}
\usepackage{mathtools}
\usepackage{color}
\usepackage{multirow}
\usepackage{enumitem}
\usepackage{hyperref}
\usepackage[T1]{fontenc}
\usepackage{scalerel,stackengine}
\usepackage[style=numeric,backend=biber]{biblatex}
\usepackage{todonotes}

\newtheorem{Lemma}{Lemma}[section]

\newtheorem{Proposition}[Lemma]{Proposition}

\newtheorem{Remark}[Lemma]{Remark}

\newcommand{\shortminus}{\scalebox{0.5}[1.0]{\( - \)}}
\DeclareMathOperator{\tr}{tr} 
\DeclareMathOperator{\cn}{cn}
\DeclareMathOperator{\sn}{sn}
\DeclareMathOperator{\dn}{dn}

\DeclareMathOperator{\const}{const}

\DeclareMathOperator{\sech}{sech}

\def\Xb{\textbf{X}}
\def\eps{\varepsilon}
\def\Ab{\textbf{A}}
\def\Bb{\textbf{B}}
\def\bb{\textbf{b}}
\def\Xib{\mathbf{\Xi}}
\def\Fib{\mathbf{\Phi}}
\def\Fb{\textbf{F}}
\def\Mb{\textbf{M}}

\graphicspath{ {./sections/pics/} } 

\usepackage{subfiles} 

\title{Time averages and periodic attractors at high Rayleigh number for Lorenz-like models}

\author[F1.IO]{Ivan Ovsyannikov}
\author[F2.JR]{Jens Rademacher}
\author[F3.RW]{Roland Welter}
\author[F4.BL]{Bing-ying Lu}
\address{Universität Hamburg\\
	Fachbereich Mathematik \\
	Bundesstraße 55 \\
	20146 Hamburg, DE}
\email{iovsyannikov@constructor.university, jens.rademacher@uni-hamburg.de, roland.welter@uni-hamburg.de, blu@sissa.it}
\thanks{\textbf{Acknowledgments:} This paper is a contribution to the project M7 of the Collaborative Research Centre TRR 181 "Energy Transfers in Atmosphere and Ocean" funded by the Deutsche Forschungsgemeinschaft (DFG, German Research Foundation) - Projektnummer 274762653. IO is also supported by the grant of the Russian Science Foundation 19-11-00280. A large part of this research was conducted while JR, BL and IO were affiliated with the University of Bremen.  We thank Camilla Nobili for helpful discussions on scaling bounds and Dmitry Turaev for discussions on homoclinic bifurcations.}

\begin{document}
	\maketitle

	\begin{abstract}
		Revisiting the Lorenz '63 equations in the regime of large of Rayleigh number, we study the occurrence of periodic solutions and quantify corresponding time averages of selected quantities. Perturbing from the integrable limit of infinite $\rho$, we provide a full proof of existence and stability of symmetric periodic orbits, which confirms previous partial results. Based on this, we expand time averages in terms of elliptic integrals with focus on the much studied average `transport', which is the mode reduced excess heat transport of the convection problem that gave rise to the Lorenz equations. We find a hysteresis loop between the periodic attractors and the non-zero equilibria of the Lorenz equations. These have been proven to maximize transport, and we show that the transport takes arbitrarily small values in the family of periodic attractors. In particular, when the non-zero equilibria are unstable, we quantify the difference between maximal and typically realized values of transport. We illustrate these results by numerical simulations and show how they transfer to various extended Lorenz models. 
	\end{abstract}	
	
	\section{Introduction}
	
	In physical processes, infinite time-averaged quantities are often of more interest than particular solutions.  Their dependence on parameters is of fundamental theoretical interest and also practical importance. The derivation of  bounds for such averages by algebraic optimization has received increasing attention in recent years \cite{FantuzziGoluskin,Goluskin2018,Goluskin2021,OlsonDoering2022}.  Bounds derived by estimates need not be sharp, but also sharp bounds may be misleading, if they are realized by dynamically unstable states \cite{wen_goluskin_doering_2022,WenDingChiniKerswell2022}. Inspired by \cite{SouzaDoering1,Goluskin2018,Goluskin2021}, in this paper we study such a situation in the regime of large Rayleigh number $\rho\gg 1$ for the famous Lorenz equations and variants thereof.  The Lorenz equations, given as follows,
	\begin{equation}\label{Lorenz63}
		\begin{array}{l}
			X' = \sigma (Y -X) \\
			Y' = \rho X - Y -XZ  \\
			Z' = -\beta Z + XY,
		\end{array}
	\end{equation}
	arose first in the context of atmospheric convection.  In the original derivation \cite{Lorenz}, Lorenz considered a fluid in a periodic box, heated from below and cooled from above, and obtained \eqref{Lorenz63} from a PDE model by retaining only the lowest order Fourier modes. The parameters in \eqref{Lorenz63} stem from the PDE, where $\sigma >0$ is the Prandtl number, characterizing the viscosity of the fluid, and where $\beta > 0$ is a shape parameter, measuring the ratio of the length of the box to its height. Often of primary interest is the parameter $\rho \geq 0$ which is the rescaled Rayleigh number, measuring the intensity of the heating.  The rescaling is chosen such that a bifurcation occurs at $\rho =1$, and indeed the Lorenz model accurately captures the onset of steady atmospheric convection rolls for $\rho >1$.  For larger values of $\rho$, the Lorenz equations do not accurately capture the full PDE model, but the relation between the hierarchy of higher order mode truncations with the convection PDE model continue to be of interest, e.g.\ \cite{Felicio_2018,Park2021,Goluskin2021,OlsonDoering2022}. 
	
	Although physically unrealistic for large regions of parameter space, \eqref{Lorenz63} is frequently used as a benchmark and testbed for nonlinear dynamics, in particular in the famous chaotic regime, but also in the context of time averages \cite{Goluskin2021}. Of particular interest is the average
	\begin{equation} \label{HeatTransportDef} 
		H(\rho,\beta, \sigma, \textbf{X}_0) = \limsup_{t \to \infty} \frac{1}{t} \int \nolimits_0^t X(t) Y(t) dt, 
	\end{equation}
	which we refer to as \emph{transport}, following \cite{SouzaDoering1}. The quantity $H$ is the mode truncated form of excess heat transport, which defines the Nusselt number by a scaling factor and constant shift. $H$ is well-defined due to the dissipation at infinity in \eqref{Lorenz63} and in particular depends on the initial condition $\textbf{X}_0$ of the solution $(X,Y,Z)(t)$ of \eqref{Lorenz63}.  The dependence of the Nusselt number on $\rho$ is of major physical interest, but is difficult to determine or bound analytically and numerically \cite{wen_goluskin_doering_2022,WenDingChiniKerswell2022}. Examining $H$ and its parameter dependence provides a tractable, non-trivial case study which can provide insight into the analysis of the Nusselt number for the full PDE. As such, an optimal bound for $H$ has been a longstanding question, which was settled by Souza and Doering \cite{SouzaDoering1}. They proved that transport is maximal in the non-trivial equilibria of \eqref{Lorenz63}, which emerge in the bifurcation at $\rho=1$. However, since these equilibria are unstable in large parameter regimes,  the question remains what values the transport takes in attractors. 
	
	The inclusion of stability in the study of transport and the resulting scaling for $\rho\gg 1$ is our main motivation for this paper.	The relation to stability is particularly clear in the case $0\leq \rho \leq 1$, where for any value of $\sigma, \beta$ the system admits a Lyapunov function and the origin is the global attractor (see for instance \cite{Sparrow}).  Hence, for any initial condition we obtain $H = 0$. 
	For higher Rayleigh numbers, the functional form of transport becomes more complicated.  At $\rho = 1$, a pitchfork bifurcation occurs, where the aforementioned non-zero fixed points $\textbf{X}_{\pm} = (\pm \sqrt{\beta (\rho - 1)}, \pm \sqrt{\beta (\rho - 1)}, \rho - 1)$ emerge.  For $\rho$ sufficiently near $1$ these non-zero fixed points seem to attract every trajectory except for $\textbf{X}_0$ belonging to the stable manifold of the origin, $W^s(\textbf{0})$, so that 
	\[ H(\rho,\beta,\sigma,\textbf{X}_0) = \left 
	\{ \begin{array}{cl} 0 & \text{if } \textbf{X}_0 \in W^s(\textbf{0}), \\ H_{\pm}(\rho,\beta):=\beta ( \rho - 1) & \text{otherwise. } \end{array} \right . \]
	Here $W^s(\textbf{0})$ is a surface of dimension 2, which means almost all initial conditions give a positive transport, namely that of the fixed points $H_{\pm}(\rho,\beta)$. This is certainly the case for initial data in their non-trivial basins of attraction. 
	
	Upon increasing $\rho$ further, additional periodic orbits emerge and further complicate the function $H$.  It has been noticed in \cite{Sparrow} that a decisive parameter for \eqref{Lorenz63} at higher $\rho$ is 
	\begin{equation*}
		\lambda = \frac{\sigma + 1}{\beta + 2}.
	\end{equation*} 
	For $0 < \lambda \leq 1$, the fixed points $\textbf{X}_{\pm}$ are locally stable for all $\rho$, cf. \cite{Sparrow} so that the fixed point transport $H_{\pm}$ is observed at least for $\textbf{X}_{0}$ belonging to a basin of attraction of positive measure.  On the other hand for $\lambda > 1$, the fixed points are locally stable only for $1 < \rho < \rho^* = \sigma(\sigma + \beta + 3)/(\sigma - \beta - 1)$.  As $\rho$ is increased through $\rho^*$, the fixed points $\textbf{X}_{\pm}$ lose stability via a sub-critical Andronov-Hopf bifurcation and at least for some open set containing $\beta = 8/3$, $\sigma = 10$, generic initial conditions give chaotic solutions \cite{Tucker}.
	
\begin{figure}
		\begin{center}
			\begin{tabular}{cc}
				\includegraphics[width=0.4\textwidth]{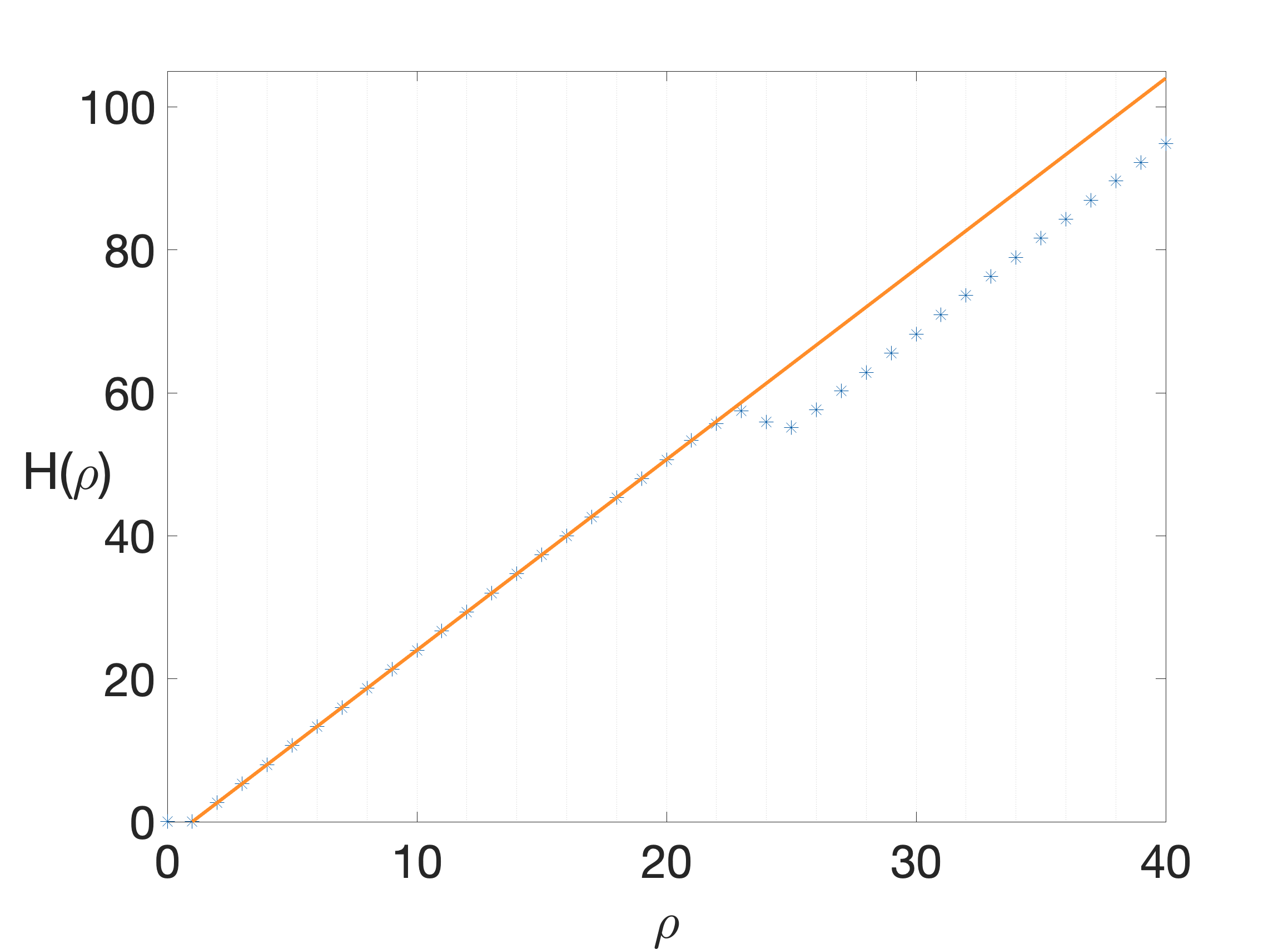}
				& \includegraphics[width=0.52\textwidth]{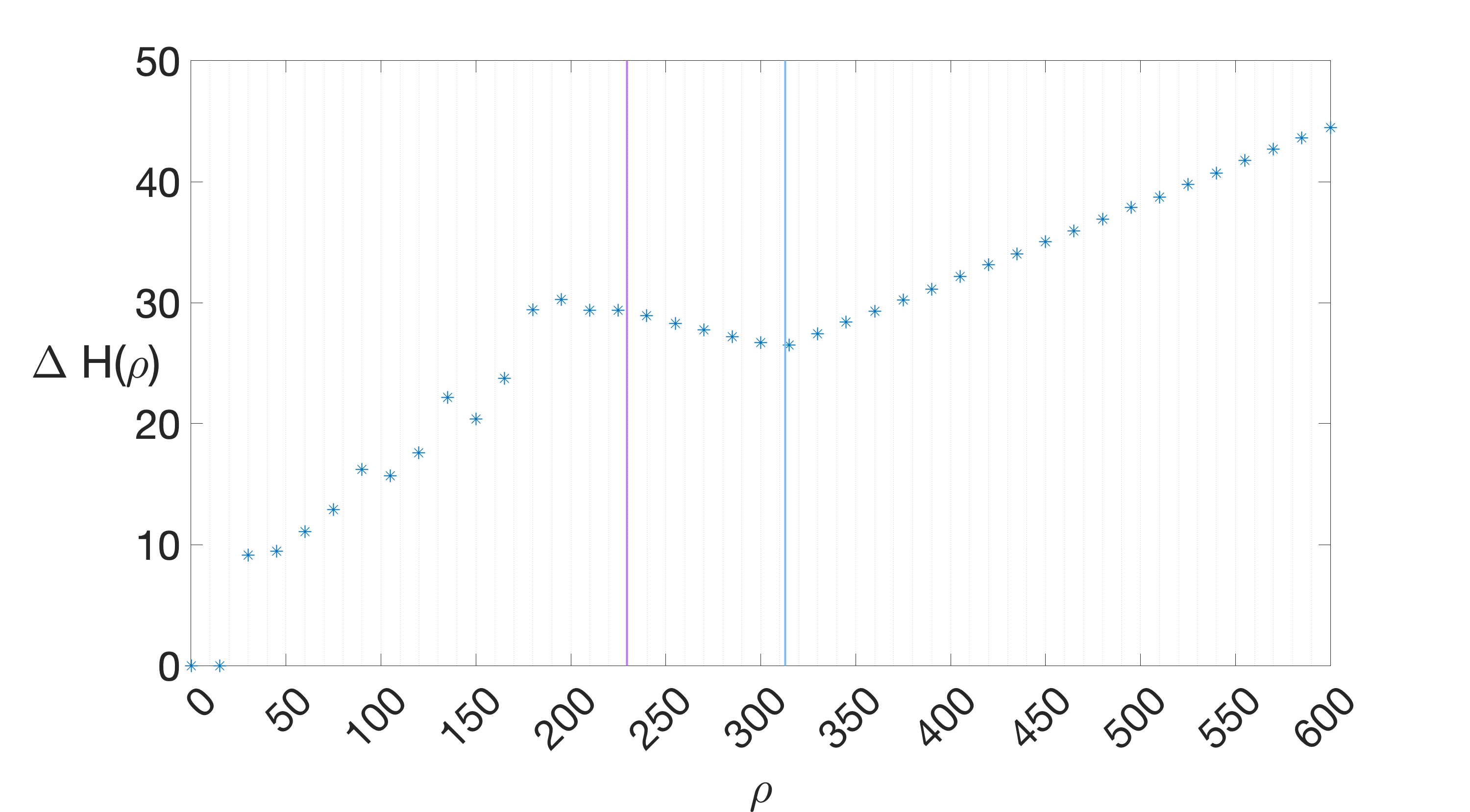} \\
				(a) & (b)
			\end{tabular}
			\caption{(a) Plot of the transport $H(\rho,\beta,\sigma,\textbf{X}_0)$ for $\beta = \frac{8}{3}$, $\sigma = 10$ near the transition to chaos $\rho \approx \rho^*$ depicting the steady convection rate $H_{\pm}(\rho,\beta)$ (orange) and the rate obtained from a numerical simulation of (\ref{Lorenz63}) with randomly selected initial conditions $\textbf{X}_0$ (blue). (b) Plot of the transport gap $\Delta H(\rho)$ for $\beta = \frac{8}{3}$, $\sigma = 10$ and for $\rho$ large. The vertical lines mark bifurcations that bound the region of `large $\rho$' as in Fig.~\ref{f:largelam}(b), in particular for $\rho$ above the blue line the transport is dominated by stable periodic orbits.}
			\label{f:SouzaDoeringBound}
		\end{center}
	\end{figure}

	As mentioned, it was proven in \cite{SouzaDoering1} that despite this complexity, for all $\rho > 1$ and any $\sigma, \beta \in \mathbb{R}_{>0}$, $\textbf{X}_0 \in \mathbb{R}^3$
	one has the simple bound 
	\[ H(\rho,\beta,\sigma,\textbf{X}_0) \leq H_{\pm}(\rho,\beta). \]
	For the Lorenz equations this bound is actually sharp as it is realized by the steady states $\Xb_\pm$. However, since $\Xb_\pm$ are unstable for $\rho>\rho^*, \lambda>1$, the transport that is realized by typical solutions might be much lower. Indeed, numerical experiments presented in \cite{SouzaDoering1} indicate that for $\rho>\rho^*$ a gap $\Delta H = H_{\pm} - H>0$ occurs,  cf.\ Fig.~\ref{f:SouzaDoeringBound}. To the best of our knowledge, quantitative results for the size of the observed transport gap for large Rayleigh number that we provide in this paper have not yet been established previously. 
	
	It has been observed already in \cite{Robbins} that for sufficiently large $\rho$ the chaotic attractor collapses and a periodic attractor occurs, but it appears that the implications for transport and other time averages have not been studied. Moreover, we found the arguments given in \cite{Robbins} and also \cite{LiZhang} to be incomplete. Sparrow \cite{Sparrow} devotes a chapter to the regime of large $\rho$ and obtains various results based on a formal application of the method of averaging. In this paper we rigorously confirm several of these results and provide a complete proof of the existence and stability of symmetric periodic attractors $\Xb_{\rm sym}=\Xb_{\rm sym}(\rho,\beta,\sigma)$ for sufficiently large $\rho$ and its dependence on $\lambda>2/3$. We also prove the existence and instability of a pair of asymmetric periodic orbits and obtain some results on homoclinic orbits.  Our analysis is based on the observation that the well known limit system as $\rho\to\infty$ possesses a Hamiltonian structure that allows one to apply the extended Melnikov theory of Wiggins and Holmes \cite{WHper}.
	
	For the periodic orbits it becomes tractable to compute the transport analytically and we quantify the gap $\Delta H$ to leading order: We show that $H_{\rm sym}(\rho,\beta,\sigma):=H(\rho,\beta,\sigma,\Xb_{\rm sym})$ is a monotone increasing function of $\lambda$, for fixed $\rho\gg1$, whose range is a subinterval of $(0,H_\pm)$ that limits to this interval as $\rho\to\infty$. Specifically, $H(\rho,\beta,\sigma,\Xb_{\rm sym})\sim\rho$  for $\rho\gg1$, just as $H_\pm$, but with a $\lambda$-dependent downshift that can bring it arbitrarily close to zero, and we provide the leading order term of the downshift in terms of elliptic integrals. 
	The resulting bifurcation diagram contains a hysteresis loop between $\Xb_\pm$ and $\Xb_{\rm sym}$ in terms of $\lambda$. This highlights difficulty to recover from low transport when $\lambda$ grows beyond the `tipping point' at $\lambda=1$.  Moreover, for any $\gamma\in(0,1]$ we can choose $\sigma(\rho)$, so that $H_{\rm sym}(\rho,\beta,\sigma(\rho))\sim \rho^\gamma$ along the periodic attractors $\Xb_{\rm sym}$.

	We employ numerical pathfollowing to corroborate the analytical results for large fixed $\rho$ and find that for large $\lambda$ the symmetric periodic orbits terminiate in a symmetric heteroclinic cycle, akin to cycles found in a different regime in \cite{Sparrow}.  We also compute the stabillity boundary of the symmetric orbits in the $(\lambda,\rho)$-plane, and find that it extends to values of $\rho$ below $200$. In Fig.~\ref{f:SouzaDoeringBound} this region begins near $\rho = 313$. It is well known that beyond this boundary various period doubling bifurcations occur \cite{Robbins}. 
	
	Finally, we turn to variants and extensions of the Lorenz equations and identify regimes in which our analytical results system remain valid and explicitly illustrate this for the Lorenz-Stenflo system. 
	
	\section{Periodic orbits at large Rayleigh number}\label{s:Melnikov}

	It is well known that \eqref{Lorenz63} possesses a semi-Lyapunov function at infinity, cf.\ \cite{SouzaDoering1}, and thus a bounded trapping region exist. However, $\Xb_\pm$ grow unboundedly with $\rho$, and as pointed out by Sparrow \cite[Chapter 7]{Sparrow}, see also \cite{Robbins}, a suitable scaling of the variables with respect to $\rho$ is given by
	\begin{equation} \label{rescaling}
		\varepsilon = \rho^{-\frac{1}{2}}, \;\; X = \varepsilon^{-1} \xi, \;\; Y = \varepsilon^{-2}\sigma^{-1} \eta, \;\; Z = \varepsilon^{-2} (\sigma^{-1} \zeta + 1), \;\;
		t = \varepsilon \tau,
	\end{equation}
	which yields the equivalent system
	\begin{equation}\label{Lorenz_resc}
		\begin{array}{l}
			\dot \xi = \eta - \varepsilon \sigma \xi \\
			\dot \eta = -\xi \zeta  - \varepsilon \eta  \\
			\dot \zeta = \xi \eta - \varepsilon \beta (\zeta + \sigma).
		\end{array}
	\end{equation}
	Notably, the symmetry $(X,Y,Z) \mapsto (-X,-Y,Z)$ of \eqref{Lorenz63} turns into $(\xi,\eta,\zeta) \mapsto (-\xi,-\eta,\zeta)$ in \eqref{Lorenz_resc}.
	
	\subsection{Hamiltonian structure and the extended Melnikov theory}
	
	It is well known that the limiting system at $\eps=0$, 
	\begin{equation}\label{Lorenz_resc_lim}
		\begin{array}{l}
			\dot \xi = \eta  \\
			\dot \eta = -\xi \zeta  \\
			\dot \zeta = \xi \eta, 
		\end{array}
	\end{equation}
	is integrable with the conserved quantities
	\begin{equation}\label{integrals}
		A = \frac{1}{2} \xi^2 - \zeta ,  \;\; B = (\eta^2 + \zeta^2)^{1/2},
	\end{equation}
	and we recap known results, essentially as provided in \cite{Sparrow}, in preparation of the existence and stability proofs.
	
	The character of a solution is completely determined by its location in the $(A,B)$-half-plane with $B\geq 0$.  In particular, $B = 0$ consists of a line of equilibria.  When $B >0$ there are two domains with distinct behavior, $D_1 =\{ 0 <|A| < B\}$, $D_2= \{A>B\}$, and two boundaries $D_3=\{ A = B\}$, $D_4= \{A = -B\}$, since the region $A < -B$ has no real solutions.  $D_4$ has the simplest solutions, consisting only of the line of equilibria $\xi = \eta = 0$, $\zeta = B$.  In other regions, the solutions are given in terms of Jacobi elliptic functions and complete elliptic integrals.  Following  \cite{ByrdFriedman}, for a given elliptic modulus $0 < k < 1$, the complete elliptic integrals of the first and second kind are defined by
	\begin{equation} \label{CompleteElliptic} K(k) = \int \nolimits _0^1 \frac{1}{\sqrt{(1-t^2)(1-k^2t^2)}} dt \hspace{,5 cm} \text{ , } \hspace{.5 cm} E(k) = \int \nolimits _0^1 \sqrt{\frac{1-k^2t^2}{1-t^2}} dt. \end{equation}
	For the Jacobi elliptic functions, one first defines the amplitude function $\mathsf{am}(u,k)$ as an inverse function via 
	\[ u = \int \nolimits_{0}^{\phi} \frac{d\theta}{\sqrt{1-k^2 \sin^2 \theta}} \hspace{,5 cm} \Leftrightarrow \hspace{.5 cm} \mathsf{am}(u,k) = \phi \]
	and the Jacobi elliptic functions are then defined
	\[ \sn (u,k) = \sin \big [ \mathsf{am} (u,k)  \big ] \hspace{.5 cm} \text{ , } \hspace{.5 cm} \cn (u,k) = \cos \big [ \mathsf{am} (u,k)  \big ] \hspace{.5 cm} \text{ , } \hspace{.5 cm} \dn (u,k) = \Big ( 1 - k^2 \sin^2 \big [ \mathsf{am} (u,k)  \big ] \Big )^{1/2}. \]
	Each $(A,B)\in D_1$ defines a symmetric periodic orbit $L_1^{A, B} = (\xi_1(\tau), \eta_1(\tau), \zeta_1(\tau))$ with period $T_1 = 4 K B^{\shortminus \frac{1}{2}}$, which can be written in terms of Jacobi elliptic functions as
	\begin{equation}\label{Jacobi_1}
		\begin{split}
			u & = \sqrt{B} \tau \text{ , }  k_1^2 = \displaystyle \frac{A + B}{2B} \\ \xi_1(\tau) & = 2 k_1 \sqrt{B} \cn (u,k_1) \\
			\eta_1(\tau) & = - 2 k_1 B \dn (u,k_1) \sn (u,k_1) \\
			\zeta_1(\tau) & = B(1 - 2 k_1^2  \sn^2 (u,k_1)).
		\end{split}
	\end{equation}
	Each $(A,B)\in D_2$ defines a pair of asymmetric periodic orbits $L_{2,\pm}^{A, B} = (\xi_{2,\pm}(\tau), \eta_{2,\pm}(\tau), \zeta_2(\tau))$ with period $\displaystyle T_2 = 4 K(k_2) k_2 B^{\shortminus \frac{1}{2}}$
	that can be represented in elliptic functions as 
	\begin{equation}\label{Jacobi_2}
		\begin{split}
			u & = \sqrt{B}k_2^{-1} \tau \text{ , }  k_2^2 = 2 \left(1+\frac{A}{B}\right)^{-1}  \\ \xi_{2,\pm}(\tau) & =  \pm 2 \sqrt{B} k_2^{-1} \dn (u,k_2) \\
			\eta_{2,\pm}(\tau) & = \mp 2 B \sn (u,k_2) \cn (u,k_2) \\
			\zeta_2(\tau) & = B(1 - 2 \sn^2 (u,k_2)).
		\end{split}
	\end{equation}
	The region $D_3: \{ A = B\}$ corresponds to the line of saddle equilibria $\xi = \eta = 0, \;\; \zeta = -B$, each with a pair of homoclinic orbits $L_{3,\pm}^{B} = (\xi_{3,\pm}(\tau), \eta_{3,\pm}(\tau), \zeta_3(\tau))$ contained in $D_3$ and given by 
	\begin{equation}\label{TanhSech}
		\begin{split}
			\displaystyle u & = \sqrt{B} \tau \\ \xi_{3,\pm}(\tau) & =  \pm 2 \sqrt{B} \sech(u) \\
			\displaystyle \eta_{3,\pm}(\tau) & = \mp 2B \tanh(u) \sech(u)  \\
			\displaystyle \zeta_3(\tau) & = B \big ( 1 - 2 \tanh^2(u) \big )
		\end{split}
	\end{equation}
	Due to the reflection symmetry of \eqref{Lorenz_resc_lim}, without loss of generality we consider only $L_{2}^{A, B} = L_{2,+}^{A, B}$ and $L_{3}^{B} = L_{3,+}^{B}$.  
	
	\medskip
	In the following we suppress the dependence of the elliptic functions and integrals on the elliptic modulus $k=k_i$, e.g., writing $\sn{u}$ for $\sn(u, k_i)$ and $K$ for $K(k_i)$. 
	
	\medskip
	Next we deviate from the approach of Sparrow, who proceeds with a formal use of the method of averaging, and instead follow that of Li and Zhang \cite{LiZhang} with corrections.  In order to exploit the Hamiltonian structure of \eqref{Lorenz_resc_lim} that is available when $\varepsilon = 0$, we introduce polar coordinates with $B$ from \eqref{integrals} given by
	\begin{equation}\label{polar_coord}
		\zeta = B \cos\varphi, \;\; \eta = B \sin\varphi, \;\; \xi = \xi
	\end{equation}
	which transform \eqref{Lorenz_resc} into
	\begin{equation}\label{Ham_1}
		\begin{array}{l}
			\dot \xi = B \sin \varphi - \varepsilon \sigma \xi \\
			\displaystyle \dot \varphi = -\xi + \varepsilon \frac{\sin \varphi}{B}((\beta - 1) B \cos \varphi + \beta \sigma)\\
			\dot B = -\varepsilon (B + (\beta - 1)B  \cos^2 \varphi + \beta \sigma \cos \varphi).
		\end{array}
	\end{equation}
	Notably, for $\varepsilon > 0$ the radial variable $B$ is no longer a conserved quantity, and must be included as a dynamical variable.  
	
	We will use that (\ref{Ham_1}) has the form
	\begin{equation}\label{Ham_2}
		\begin{array}{l}
			\dot \xi = f_1(\xi, \varphi, B) + \varepsilon g_1 (\xi, \varphi, B) \\
			\dot \varphi = f_2(\xi, \varphi, B) + \varepsilon g_2 (\xi, \varphi, B) \\
			\dot B = \varepsilon g_3 (\xi, \varphi, B).
		\end{array}
	\end{equation}
	for smooth functions $f_i$, $g_i$. In this formulation, at $\varepsilon = 0$, the first two equations possess the Hamiltonian structure with $A(\xi,  \varphi, B)$ defined in (\ref{integrals}) serving as a Hamiltonian:
	\begin{equation}\label{Ham_func}
		A(\xi,  \varphi, B) = \displaystyle \frac{1}{2}\xi^2 - B \cos\varphi, \;\; \displaystyle f_1 = \frac{\partial A}{\partial \varphi}, \;\;  f_2 = -\frac{\partial A}{\partial \xi}.
	\end{equation} 
	As noted above, the solutions at $\varepsilon = 0$ are given by families of periodic orbits, saddle equilibria and their homoclinic orbits. Therefore, as already noticed in \cite{LiZhang}, we can apply the extended Melnikov perturbation theory of Wiggins and Holmes  \cite{WHper, WHhom} in order to identify which periodic and homoclinic orbits of (\ref{Ham_1}) persist for small $\varepsilon > 0$. The main idea of this method is that for $\varepsilon = 0$ the phase space is represented as a one-parametric family of two-dimensional manifolds (parametrized by $B$ in our case), and on each manifold the system is Hamiltonian. Then, generically, in the phase space there exist two-parameter families of periodic orbits, parametrized by $A$ (the Hamiltonian) and $B$, and one-parameter families of homoclinic orbits. Upon perturbing by $\eps>0$, this structure is destroyed, and, generically, one can expect that only isolated periodic orbits and isolated homoclinic orbits will exist. The existence, local uniqueness and the topological type of these objects can be established with the help of Melnikov integrals. 
	
	For the periodic orbits $L^{A,B}_{i}$, the analysis simplifies when changing canonical Hamiltonian variables  for $\varepsilon = 0$  from $(\xi, \varphi)$ to action-angle variables $(I, \theta)$. The action variable $I$ is computed as 
	\begin{equation}\label{action_sym}
		I(A, B) = \oint_{L_i^{A, B}} \varphi d \xi = -\oint_{L_i^{A, B}} \xi d \varphi = \int\limits_0^{T_i} \xi^2 d\tau,
	\end{equation}
	and the angle $\theta\in[0,1)$ increases from an $0$ to $1$ along the periodic orbits with constant frequency $\Omega(A, B)$ given by
	\begin{equation}\label{frequency_sym}
		\frac{1}{\Omega(A, B)} = \frac{\partial I}{\partial A}\Bigr|_B.  
	\end{equation}
	Here and below we use the notation $\displaystyle \frac{\partial R}{\partial P}\Bigr|_Q$ from thermodynamics and elsewhere to emphasize that we differentiate the quantity $R$ only with respect to its explicit dependence on $P$, neglecting implicit relationships between $P$ and $Q$.\\ 
	In these new variables (\ref{Ham_1}) turns into 
	\begin{equation}\label{action-angle}
		\begin{array}{l}
			\dot I = \varepsilon F(I, \theta, B)  \\
			\dot \theta = \Omega + \varepsilon G(I, \theta, B) \\
			\dot B = \varepsilon \tilde g_3(I, \theta, B),
		\end{array}
	\end{equation}
	where
	\begin{equation}\label{act-ang_func}
		\begin{array}{l}
			\displaystyle F(I, \theta, B) = \frac{\partial I}{\partial \xi} \tilde g_1 +
			\frac{\partial I}{\partial \varphi} \tilde g_2 + \frac{\partial I}{\partial B} \tilde g_3 \\[2ex]
			\displaystyle G(I, \theta, B) = \frac{\partial \theta}{\partial \xi} \tilde g_1 +
			\frac{\partial \theta}{\partial \varphi} \tilde g_2 + \frac{\partial \theta}{\partial B} \tilde g_3 \\[2ex]
			\tilde g_i(I, \theta, B) = g_i(\xi(I, \theta, B), \varphi(I, \theta, B), B), \;\; i = 1, 2, 3.
		\end{array}
	\end{equation}
	For small $\varepsilon > 0$, a trajectory starting from the initial point $(I_0, 0, B_0)$ has monotonically increasing  $\theta$ from $0$ to $1$ with velocity $\varepsilon$-close to $\Omega$, and slowly evolving $(I, B)$ with velocity $O(\varepsilon)$. Hence, the trajectory stays in an $\varepsilon$-neighborhood of the corresponding unperturbed periodic trajectory, and reaches $\theta = 1$ after finite time $T_\varepsilon = 1/\Omega + O(\varepsilon)$, thus returning to the starting plane $\{\theta=0\}$. This defines a Poincar\'e map $(I_0, B_0) \to (I_1, B_1)$, 
	
	\begin{equation}\label{Poincare}
		(I_1, B_1) = (I_0, B_0) + \varepsilon (M_1, M_3) + O(\varepsilon^2),
	\end{equation}
	whose fixed points correspond to periodic orbits,  and where $M_1, M_3$ will be Melnikov integral terms. We denote the linearization matrix at $(I_0, B_0)$ as
	\[  DM(I_0, B_0) =  \displaystyle \frac{\partial(M_1, M_3)}{\partial(I, B)}(I_0, B_0). \]
	Theorem 3.2 of \cite{WHper} states that for any $(I_0,B_0)$ for which $M_1 = M_3 = 0$, $\det DM(I_0,B_0) \neq 0 $ there exists $\eps_0>0$ and an isolated fixed point of the Poincar\'e map \eqref{Poincare} in an $\varepsilon$-neighbourhood of $(I_0,B_0)$ for any $0<\eps<\eps_0$. This corresponds to a persistent periodic orbit of \eqref{Lorenz_resc} and \eqref{Lorenz63} for $0<\eps\ll1$. Moreover, in case $(M_1,M_3)\neq 0$ for $0<\eps\ll1$ there is no periodic orbit in a neighbourhood of $(I_0,B_0)$.  The expressions for $M_1$ and $M_3$ are given in \cite{WHper}, and in our notation read
	\begin{equation}\label{Melnikov_per}
		\begin{array}{l}
			\displaystyle M_1 = \frac{1}{\Omega} \tilde M_1 + 
			\frac{\partial I}{\partial B}\Bigr|_A M_3, \;\; \;
			\tilde M_1 = \int\limits_0^{T_i} \left[ f_1 \tilde g_2 - f_2 \tilde g_1 + \frac{\partial A}{\partial B}\Bigr|_{\xi, \varphi} \tilde{g}_3 \right] (L_i^{A,B}(\tau)) d\tau \\
			\displaystyle M_3 = \int\limits_0^{T_i} \tilde g_3(L_i^{A,B}(\tau)) d \tau; 
		\end{array}
	\end{equation}
	in particular $(M_1,M_3)=0$ is equivalent to $(\tilde M_1, M_3) = 0$. We emphasize that we apply \cite[Theorem 3.2]{WHper} separately for two topologically different families of periodic orbits, lying respectively in domains $D_1$ and $D_2$. For every fixed point of the Poincar\'e map $(I_0, B_0) \in D_i$ there exists small $\eps_0>0$ such that for $0 < \varepsilon < \varepsilon_0$ the corresponding periodic orbit lies entirely in $D_i$.
	Also, it is clear that upon increasing $\varepsilon$ further, the trajectory can touch the boundary of both domains -- $D_3$, and cross it. This scenario 
	is confirmed below by numerical experiments in section~\ref{sec:numerics}, see Fig.~\ref{f:ABproj}.
	
	The stability type of such a fixed point, and thus the periodic orbit, is determined by the eigenvalues $\nu_\pm$ of $DM(I_0,B_0)$ are given by
	\begin{equation}\label{characteristic}
		\nu_\pm = 1 + \frac{\varepsilon}{2} \big ( \tr DM \pm \sqrt{(\tr DM)^2 - 4 \det DM} \big ) + O(\varepsilon^2).
	\end{equation}
	Although $(M_1, M_3)$ depend on $(I,B)$, it is more convenient to consider them as functions of $k_i, B$, where $k_i$ are the moduli defined in (\ref{Jacobi_1}), (\ref{Jacobi_2}).  The following simplified formulas for the trace and determinant can be obtained by using the fact that the Melnikov functions are zero at the appropriate value of $(I_0,B_0)$, and by changing variables:
	\begin{equation}
		\label{derivatives}
		\begin{split}
			\tr DM & = \left( \frac{\partial I}{\partial k_i}\Bigr|_{B}\right)^{-1} \Big [ \frac{1}{\Omega} \frac{\partial \tilde M_1 }{\partial k_i} \Bigr|_{B} +
			\frac{\partial I}{\partial B}\Bigr|_A  \frac{\partial M_3 }{\partial k_i} \Bigr|_{B} \Big ] + \frac{\partial M_{3}}{\partial B}\Bigr|_{k_i} - \frac{\partial M_{3}}{\partial k_i}\Bigr|_{B} \frac{\partial I}{\partial B}\Bigr|_{k_i} \left( \frac{\partial I}{\partial k_i}\Bigr|_{B}\right)^{-1} \\
			\det DM & = \frac{1}{\Omega} \left( \frac{\partial I}{\partial k_i}\Bigr|_{B}\right)^{-1} \Big [ \frac{\partial \tilde{M}_1}{\partial k_i} \Bigr|_{B} \frac{\partial M_{3}}{\partial B}\Bigr|_{k_i} - \frac{\partial \tilde{M}_1}{\partial B} \Bigr|_{k_i} \frac{\partial M_{3}}{\partial k_i}\Bigr|_{B}  \Big ].
		\end{split}
	\end{equation}
	\begin{Remark} 
		We briefly comment on the previous existence and stability studies in the literature. 
		The authors of \cite{LiZhang} also follow the method from \cite{WHper} with the difference that in formula (\ref{polar_coord}) they define the radius as $(B + \rho)$, where $B = \const$ and $\rho$ being the third dynamical variable in system (\ref{Ham_1}) (not the Rayleigh number, as in our notation).  However, in the new coordinates they then write formula (\ref{Ham_func}) as
		$$
		A(\xi,  \theta, \rho) = \displaystyle \frac{1}{2}\xi^2 - B\cos\theta,
		$$
		which is incorrect as the term $\displaystyle \frac{\partial A}{\partial \rho}$ is missing from the Melnikov integrals (\ref{Melnikov_per}). This is why the results of \cite{LiZhang} differ from   \cite{Sparrow}, \cite{Robbins} and our results.
		
		In \cite{Robbins} the Lorenz system of another form is considered. It can be obtained from system (\ref{Lorenz63}) via a coordinate transformation and setting $\beta=1$, thus the parameter space is reduced. Moreover, there is a gap in the argumentation. Namely, for the unperturbed periodic orbit $(x_0(t), y_0(t), z_0(t))$ a perturbation $(x_1(t), y_1(t), z_1(t)) = O(\varepsilon)$ is considered. The author solves the system of differential equations for $(x_1(t), y_1(t), z_1(t))$ under an assumption $z_0(t) \neq 0$ (see  \cite[formula (6)]{Robbins}), but on the symmetric orbit function $z_0(t)$ obviously vanishes two times. Thus, while the results of the computation seem correct, they are not sufficiently justified.
		
		In the book \cite{Sparrow}  system (\ref{Lorenz_resc}) is analyzed via formally averaging over the unperturbed periodic orbits. Isolated equilibrium points of the averaged system correspond then to isolated periodic orbits of the original system. Formulas (\ref{Meln_sym_3}) and  (\ref{Meln_nonsym_3}), that determine the existence and uniqueness of periodic orbits in domains $D_1$ and $D_2$, were also obtained there, however without a rigorous proof of monotonicity. Also, analogues of (\ref{tr_stable}) and (\ref{tr_unstable}) were obtained in \cite{Sparrow}, giving the sign of the trace of the linearization matrix. However, the determinant was not computed, and thus the stability of the symmetric periodic orbit and instability of the pair of non-symmetric periodic orbits could not be determined.
	\end{Remark}
	\begin{Remark}\label{rem:hom}
		The question of persistence also arises for the homoclinic orbits \eqref{TanhSech} of \eqref{Ham_1} in region $D_3$. From the line of the corresponding saddle equilibria, for $0<\varepsilon \ll 1$ only one saddle equilibrium $(\xi, \varphi, B) = (0, \pi, \sigma)$ remains. Hence, it is possible that its stable and unstable manifolds will form homoclinic orbits, and periodic orbits may appear in bifurcations from the homoclinic orbits at $\eps>0$. 
		In  \cite{WHhom} such a phenomenon was studied, however, the results are not valid in the claimed generality and do not apply in our case as discussed in \S\ref{sec:hom} below.
	\end{Remark}

	\subsection{Symmetric periodic orbits ($D_1$)}\label{s:perD1}
	
	In the regime $D_1$ the explicit solutions (\ref{Jacobi_1}) have the following form in polar coordinates:
	\begin{equation*} \begin{split} 
			\xi & = 2 k_1 \sqrt{B} \mathsf{cn} (\sqrt{B} \tau) \\
			\sin \varphi & = - 2k_1 \mathsf{dn}(\sqrt{B}\tau ) \mathsf{sn} (\sqrt{B}\tau ) \\
			\cos \varphi & = 1-2k_1^2 \mathsf{sn}^2(\sqrt{B}\tau ).
	\end{split} \end{equation*}
	Hence, from (\ref{action_sym}) and (\ref{frequency_sym}) the period and the action are 
	\begin{equation}\label{e:TI_D1}
		\displaystyle
		T_1 = \frac{4K}{\sqrt{B}} = \frac{1}{\Omega(k_1, B)}   \hspace{.5 cm} \text{ , } \hspace{.5 cm} I(k_1, B) =  16 \sqrt{B}(E - (1 - k_1^2)K).
	\end{equation}
	Substituting these explicit solutions into (\ref{Melnikov_per}) using the formulas for $f_i,g_i$ from (\ref{Ham_2}) allows to compute the Melnikov integrals explicitly as
	\begin{equation*}\begin{split}
			\tilde{M}_1 & = \int_0^{4 K} \big [ - 4 \sigma k_1^2 \sqrt{B} \mathsf{cn}^2(u ) + \sqrt{B} \beta (1 - 2 k_1^2 \mathsf{sn}^2(u)) + \frac{\beta \sigma}{\sqrt{B}} \big ] du\\ 
			&= (B\beta + \beta \sigma ) \frac{4K}{\sqrt{B}} - 16 \sigma \sqrt{B} ( E - (1-k_1^2 ) K) + 8 \beta \sqrt{B} ( E - K) \\
			M_3 & = -(B + \sigma ) \beta \frac{4K}{\sqrt{B}} + \int_0^{4 K} \big [ (4 k_1^2 (\beta-1)\sqrt{B} + \frac{2 \sigma \beta k_1^2}{\sqrt{B}}) \mathsf{sn}^2(u) - 4k_1^4 (\beta-1)\sqrt{B} \mathsf{sn}^4(u) \big ] du \\ & = -(B + \sigma ) \beta \frac{4K}{\sqrt{B}} - 8(2(\beta-1)\sqrt{B} + \frac{\sigma \beta}{\sqrt{B}}) (E-K) - \frac{16 (\beta-1)\sqrt{B}}{3} ( -2(1+k_1^2)E + (2+k_1^2) K ).
	\end{split} \end{equation*}
	Equating these expressions to zero, equivalently gives the equations 
	\begin{equation} \label{Meln_sym_2} \begin{split}
			K \beta \sigma + \beta B (2E - K ) - 4 B \sigma ( E - (1-k_1^2)K) & = 0 \\
			4B \big [ (1-k_1^2) K + (2k_1^2-1) E \big ] + 3 \beta \sigma (2 E - K) + \beta B \big [ (4k_1^2-1) K + 4(1-2k_1^2)E \big ] & = 0
	\end{split} \end{equation}
	for $A, B$ as necessary conditions for the existence of persistent periodic orbits.  Indeed, the same equations were obtained by C.\ Sparrow and P.\ Swinnerton-Dyer (see \cite[Appendix K]{Sparrow}) by formally using the method of averaging. In fact, part of the following analysis is an extension of that in \cite{Sparrow} regarding \eqref{Meln_sym_2}.
	
	We are now ready to formulate and prove our main result for symmetric periodic orbits. 
	
	\begin{Proposition}
		For every $\lambda = \frac{\sigma + 1}{\beta + 2} > 2/3$ and $\varepsilon$ sufficiently small there exists a (locally) asymptotically stable symmetric periodic orbit. It is the unique periodic orbit in an order $\eps$-neighborhood of an explicit solution $L_1^{A, B}$, where $(A,B)$ is the unique solution of  \eqref{Meln_sym_2} for the chosen value of $\lambda$.
	\end{Proposition}
	\begin{proof}
		We start by showing that (\ref{Meln_sym_2}) possesses a unique solution in the terms of $k_1$ and $B$.  In order to solve \eqref{Meln_sym_2}, we first show that the coefficient of $B$ in the second equation is strictly positive for $0 \leq k_1 \leq 1$.  Let $e_1,e_2$ be defined by 
		\begin{equation} \label{EllipticIntegralExp1} e_1(k_1) = (1-k_1^2) K + (2k_1^2-1) E \hspace{1 cm} \text{ , } \hspace{1 cm} e_2(k_1) = (4k_1^2-1) K + 4(1-2k_1^2)E  \end{equation}
		For $e_1$ write
		\[ e_1(k_1) = (1-k_1^2) K + (2k_1^2-1) E = (K - E) (1-k_1^2) + E k_1^2 \]
		and since it is clear from the definitions (\ref{CompleteElliptic}) that $K \geq E \geq 0$, this quantity is strictly positive except when $k_1 = 0$, in which case it is zero.  The function $e_2$ is also strictly positive, although it is more work to prove it.  Rather than interrupt the Melnikov analysis, we include this proof in Appendix \ref{app:EllipticIntegralPositivity}.  Thus, in order to have positive solutions for $B$, one should have $2E - K < 0$.  As noticed in \cite[Appendix K]{Sparrow}, this  is possible precisely when $k_* < k_1 \leq 1$ with $k_* \approx 0.908909$.
		
		Since the coefficient of $B$ is strictly positive, one can solve the second equation for $B$:
		\begin{equation} \label{Meln_sym_B} B = \frac{3 \beta \sigma (K-2E)}{4 e_1 + \beta e_2 } =: \frac{3 \beta \sigma (K-2E)}{ d_1 }, \end{equation}
		where $d_1 = 4 e_1 + \beta e_2 $ denotes the denominator of the first fraction.  Substituting this into the first equation gives 
		\[ Kd_1 - 3\beta(K-2E)^2 -12 \sigma (K-2E)(E-(1-k_1^2)K) = 0 \]
		as the remaining equation. 
		Moving the term involving $\sigma$ to the right hand side, using the definition $\lambda = (\sigma+1)/(\beta+2)$, and after some arithmetic one obtains
		\begin{equation}\label{Meln_sym_3}
			\displaystyle
			2 \lambda - 1 = \frac{K((1 - k_1^2)K +  (2 k_1^2 - 1)E)}{3 (K - 2E)(E - (1 - k_1^2)K)}.
		\end{equation}
		Next, we prove the claim of \cite{Sparrow} that the right-hand side in equation (\ref{Meln_sym_3}) is a monotonically decreasing function of $k_1$ in the domain $k_* < k_1 \leq 1$ and the right hand side limits to $1/3$ at $k_1=1$. This immediately implies the existence and uniqueness of the solution of equations (\ref{Meln_sym_2}) for $\lambda > 2/3$.  Writing the right hand side of \eqref{Meln_sym_3} as 
		$$
		\displaystyle \frac{1}{3(E - (1 - k_1^2)K)}\left( K (1 - k_1^2) + \frac{E}{1 - 2 E/K} \right),
		$$ 
		we claim that the first factor and both summands in the parentheses are non-negative and monotonically decreasing. Indeed, we compute for the prefactor and the first summand that
		\begin{align*}
			&\displaystyle \frac{d}{dk_1} (E - (1 - k_1^2)K) = k_1 K > 0, &&
			\left. (E - (1 - k_1^2)K) \right|_{k_1 = k_*} = (2 k_*^2 - 1) E > 0,\\
			&\displaystyle \frac{d}{dk_1}( K (1 - k_1^2)) = \frac{E - K}{k_1} - k_1 K < 0, &&
			\lim\limits_{k_1 \to 1}(K (1 - k_1^2))  = 0,
		\end{align*}
		cf.\ \cite{ByrdFriedman}. The second summand is decreasing due the monotonicity of $E$ and $K$, and is positive since $K>2E$. Together with $E(1)=0$ and $K\to\infty$ as $k_1\to1$ the right-hand side of equation (\ref{Meln_sym_2}) is monotonically decreasing from $\infty$ to $1/3$ when $k_1$ increases from $k_*$ to $1$. This range corresponds to the values $\lambda > 2/3$ on the left-hand side of \eqref{Meln_sym_3}. Also, every $k_1 \in (k_*, 1)$ gives a unique positive value for $B$ via formula (\ref{Meln_sym_B}).
		
		Having found the unique solution to the condition $M_1 = M_3 =0$ for each $\lambda>2/3$, we now consider the determinant from \eqref{derivatives}.  Computing the derivatives in  (\ref{derivatives}) gives
		\begin{equation*}\begin{split} \det DM & = \frac{1}{3 k_1^2 B^{3}(1-k_1^2)} \Big [ \tilde{c}_2 K^2 + \tilde{c}_1 KE + \tilde{c}_0 E^2,  \Big ]   \end{split} \end{equation*}
		where $\tilde{c}_0,\tilde{c}_1,\tilde{c}_2$ are coefficients depending on $\sigma, \beta, B$ only.  Using the identity (\ref{Meln_sym_B}) to eliminate $B$ and the identity $\sigma = \lambda (\beta+2)-1$ to eliminate $\sigma$, one obtains
		\begin{equation*}\begin{split} \tilde{c}_2 K^2 + \tilde{c}_1 KE + \tilde{c}_0 E^2 & = \frac{48 \sigma^2 \beta^2 (\beta+2)}{d_1^2 } \Big [ \hat{c}_4 K^4 + \hat{c}_3 K^3E + \hat{c}_1 KE^3  + \hat{c}_0 E^4 \Big ]   \end{split} \end{equation*}
		for coefficients $\hat{c}_0,\hat{c}_1,\hat{c}_2,\hat{c}_3,\hat{c}_4$ depending only on $\lambda, \beta$.  Finally, applying (\ref{Meln_sym_3}) to eliminate $\lambda$ yields 
		\[ \hat{c}_4 K^4 + \hat{c}_3 K^3E + \hat{c}_1 KE^3  + \hat{c}_0 E^4 = \frac{d_1}{E-(1-k_1^2)K} F_1  \]
		where 
		\begin{equation*}\begin{split} F_1 = K^3 (1-k_1^2) \left(\frac{7}{20} E - (1-k_1^2) K\right) + E^3 \left ( (2k_1^2-1) E - (1-k_1^2) K \right ) + KE(1-k_1^2) \left ( (\frac{73}{20}-2k_1^2)K^2 - 6 KE +5E^2 \right ). \end{split} \end{equation*}
		This equality expresses $F_1$ as a sum of strictly positive terms, since for the first two terms we have
		\begin{align*}
			\frac{d}{dk_1} \left(\frac{7}{20} E - (1 - k_1^2)K\right) = \left(k_1 + \frac{7}{20 k_1}\right) K - \frac{7}{20 k_1} E > 0,&\quad
			\left(\frac{7}{20} E - (1 - k_1^2)K\right) \big |_{k_1 = k_*} = \left(2 k_*^2 - \frac{33}{20}\right) E > 0 \\ 
			(2k_1^2-1)E - (1 - k_1^2)K  &\geq \frac{7}{20} E - (1 - k_1^2)K > 0, \end{align*}
		and we prove positivity of the last term in Appendix \ref{app:EllipticIntegralPositivity2}.  This proves positivity of the determinant and thus existence and local uniqueness for $0<\eps\ll 1$ by \cite[Theorem 3.2]{WHper} as mentioned above.
		
		Having proven the existence of a persistent periodic orbit, we can now prove its stability.  Since we already have $\det DM > 0$ we study the trace from \eqref{derivatives}.  In region $D_1$ one has $\displaystyle \frac{\partial I}{\partial B} \big |_A = \frac{4}{\sqrt{B}}(2 E - K)$, so that by computing the derivatives in the formula (\ref{derivatives}) we find 
		\[ \tr DM = \frac{4}{3 B^{3/2}} \Big [ \hat{c}_1 K + \hat{c}_0 E \Big ], \]
		where $\hat{c}_0, \hat{c}_1$ are coefficients depending only on $\sigma$, $\beta$, $B$ and $k_1$ (but not $E$ or $K$).  Using the first equation of (\ref{Meln_sym_2}) and $\lambda>2/3$ (so that $\beta<2\sigma$) we express $E$ as
		\begin{equation}\label{E_for_tr}
			\displaystyle
			E = \frac{K ( \beta (B - \sigma) - 4 B (1 - k_1^2) \sigma)}{2 B (\beta - 2 \sigma)}
		\end{equation}
		and substituting this into the equation for the trace, we find
		\begin{equation}\label{tr_stable}
			\displaystyle \tr DM = -\frac{4 K(1 + \beta + \sigma)}{\sqrt{B}} < 0.
		\end{equation}
		Since $\tr DM < 0$ and $\det DM > 0$, from (\ref{characteristic}) it follows that for $\eps>0$ sufficiently small the eigenvalues lie inside the unit circle so that the periodic orbit is asymptotically stable.
	\end{proof}
	
	\subsection{Asymmetric periodic orbits ($D_2$)}\label{s:D2per}
	
	In the regime $D_2$ we consider the explicit solutions $L_2^{A,B}(\tau)$ given in (\ref{Jacobi_2}). In polar coordinates these are given as
	\begin{equation*}\begin{split}
			\xi & = 2 \sqrt{B}k_2^{-1} \mathsf{dn} (\sqrt{B}k_2^{-1} \tau) \\
			\sin \varphi & = - 2 \mathsf{sn}( \sqrt{B}k_2^{-1} \tau ) \mathsf{cn} (\sqrt{B}k_2^{-1} \tau ) \\
			\cos \varphi & = 1-2 \mathsf{sn}^2(\sqrt{B}k_2^{-1}\tau ) ,
	\end{split} \end{equation*}
	with period and action 
	\begin{equation}
		\displaystyle
		T_2 = \frac{4K k_2}{\sqrt{B}} = \frac{1}{\Omega(k_2, B)} ,\quad I(k_2, B) =  16 \sqrt{B} k_2^{-1} E.
	\end{equation}
	Substituting these explicit solutions into (\ref{Ham_2}) and (\ref{Melnikov_per}) we compute the Melnikov integrals explicitly as
	\begin{equation*}\begin{split}
			\tilde{M}_1 & = \int_0^{4K} \frac{k_2}{\sqrt{B}} \big [ - 4 \sigma B k_2^{-2} \mathsf{dn}^2(u) + \beta B \big ( 1-2 \mathsf{sn}^2(u) \big ) + \beta \sigma \big ] du = \frac{4 k_2}{\sqrt{B}} \big [ \beta \sigma K -  4\sigma B k_2^{-2} E + \beta B \big ( K - 2k_2^{-2}(K-E) \big ) \big ]  \\
			M_3 & = -\int_0^{4 K} \frac{k_2}{\sqrt{B}} \big [ \beta (B+\sigma) - 2 (\beta \sigma + 2 (\beta - 1)B ) \mathsf{sn}^2(u) + 4 (\beta - 1) B \mathsf{sn}^4(u) \big ] du \\ & = -\frac{\beta k_2 (B+\sigma)}{\sqrt{B}}4K + \frac{8}{\sqrt{B}}(\beta \sigma + 2 (\beta - 1)B ) \frac{K-E}{k_2} - 16 (\beta - 1) \sqrt{B} \frac{(2+k_2^2) K - 2(1+k_2^2)E}{3k_2^3}.
	\end{split} \end{equation*}
	Equating these expressions to zero, we see that we must solve
	\begin{equation}\label{Meln_nonsym_2}
		\begin{split}
			\beta \sigma K - 4 B \sigma k_2^{-2} E - \beta B k_2^{-2}((2 - k_2^2)K - 2 E) & = 0 \\
			4 B \left[(2 - k_2^2) E - 2 (1 - k_2^2) K\right] + 
			\beta B \left[4 (k_2^2 - 2) E + (3 k_2^4 - 8 k_2^2 + 8) K\right] + 3 \beta \sigma k_2^2 \left[2 E - K (2 - k_2^2)\right] & = 0.
		\end{split}
	\end{equation}
	for $A, B$.  As above, the same conditions were obtained in  \cite[Appendix K]{Sparrow} by formal averaging.  We can now state and prove our existence and instability result for asymmetric periodic orbits:
	
	\begin{Proposition}\label{p:per}
		For every $\varepsilon$ sufficiently small and all $2/3 < \lambda < 1 $ there exists a pair of saddle-type asymmetric periodic orbits. Each is unique in an order $\eps$-neighborhood of $L_2^{A,B}$, where $(A,B)$ corresponds to a unique solution $(k_2,B)$ of \eqref{Meln_nonsym_2} for the chosen value of $\lambda$. 
	\end{Proposition}
	
	\begin{proof}
		The coefficient of $B$ in the first equation of \eqref{Meln_nonsym_2} is always positive, since
		\[ (2-k_2^2) K - 2E = \frac{1}{4}\int_0^{4K} \mathsf{sn}^2(u) \mathsf{cd}^2(u) du > 0, \]
		so that the unique solution in the terms of $B$ is
		\begin{equation} \label{Meln_nonsym_B} B = \frac{\beta \sigma K k_2^2}{4 \sigma E + \beta ( ( 2-k_2^2)K-2E)} =: \frac{\beta \sigma K k_2^2}{d_2}, \end{equation}
		where $d_2$ denotes the denominator. 
		Substitution into the second equation of \eqref{Meln_nonsym_2} gives an equation for $k_2$:
		\begin{equation}
			\label{Meln_nonsym_3}
			2 \lambda - 1 = \frac{K((2 - k_2^2)E - 2 (1 - k_2^2)K)}{3 E ((2 - k_2^2)K - 2 E)}.
		\end{equation}
		Next, we prove that the right-hand side of \eqref{Meln_nonsym_3} monotonically decreases from $1$ to $1/3$ as $k_2$ increases from $0$ to $1$, which implies precisely for every $\lambda\in (2/3,1)$ there exists a unique solution $k_2$.  To prove monotonicity, first we compute
		\[ \frac{d}{d k_2} \left [ \frac{K((2 - k_2^2)E - 2 (1 - k_2^2)K)}{3 E ((2 - k_2^2)K - 2 E)}  \right ] 
		= \frac{2 P}{3 E^2 k_2 ( 1-k_2^2)((2-k_2^2)K-2E)^2 } \]
		where
		\begin{equation}          \label{NegativePolynomial}  
			\begin{split} P = & E^3 \big ( -(2-k_2^2)E + 2(1-k_2^2)K \big ) + 3 K E^2 (1-k_2^2) \big ( 2E - (2-k_2^2)K \big ) \\ 
				& + \frac{1}{2} K^2 (1-k_2^2)(2-k_2^2) \big ( 3E ( -2E + (2-k_2^2)K) + K ((2-k_2^2)E-2K ) \big ). 
		\end{split}  \end{equation}
		The three summands of $P$ are each negative since for the first we have
		\[ \frac{d}{dk_2} \big [ -(2-k_2^2)E+2(1-k_2^2)K \big ] = 2 k_2( E - K) < 0 \hspace{1 cm} \text{ and } \hspace{1 cm} \lim_{k_2 \to 0 } -(2-k_2^2)E+2(1-k_2^2)K = 0, \]
		for the second $K\geq 0$, $2E-K<0$, and for the third we estimate its nontrivial factor from above by the negative $-4K^2+2E^2$ using the quadratic estimate $2KE\leq K^2+E^2$. 
		
		Regarding the determinant, by computing the derivatives in (\ref{derivatives}) we find it has the form
		\begin{equation*}\begin{split} \det DM & = \frac{1}{3 k_2^2 B^{3}(1-k_2^2)} \Big [ \tilde{c}_2 K^2 + \tilde{c}_1 KE + \tilde{c}_0 E^2  \Big ],   \end{split} \end{equation*}
		where $\tilde{c}_0,\tilde{c}_1,\tilde{c}_2$ are coefficients depending only on $\sigma, \beta, B, k_2$.  Again using the identities (\ref{Meln_nonsym_B}), $\sigma = \lambda (\beta+2)-1$ and (\ref{Meln_nonsym_3}) to eliminate $B, \sigma, \lambda$ one finds 
		\begin{equation} \label{D2_DetIdent} \det DM = \frac{ -16 (\beta + 2) d_2 F_2 }{3B k_2^2 (1-k_2^2) KE ((2-k_2^2)K - 2E) },  \end{equation}
		where $F_2$ is given by
		\begin{equation*}\begin{split} F_2(k_2) = (2-k_2^2)E^4 - 8(1-k_2^2)E^3 K + 6(1-k_2^2)(2-k_2^2)E^2K^2 - 2(2-k_2^2)^2(1-k_2^2)E K^3 + (2-k_2^2)(1-k_2^2)^2 K^4. \end{split} \end{equation*}
		In Appendix \ref{app:EllipticIntegralPositivity3} we prove that $F_2(k_2)>0$ for $0 < k_2 < 1$, hence it follows that $\det DM < 0$.  This proves existence and local uniqueness for $0<\eps\ll 1$ by \cite[Theorem 3.2]{WHper} as mentioned above.
		
		Turning now to the question of stability, in region $D_2$ one has $\displaystyle \frac{\partial I}{\partial B} \big |_A = 4(2 E - K(2-k_2^2))/(k_2\sqrt{B})$, and hence by computing the derivatives in  (\ref{derivatives}) one obtains 
		\[ \tr DM = \frac{1}{3 B^{3/2}k_2^3} \Big [ \check{c}_1 K + \check{c}_0 E \Big ], \]
		where $\check{c}_1, \check{c}_0$ are coefficients depending only on $\sigma$, $\beta$, $B$ and $k_2$ (but not $E$ or $K$).  Using the identity (\ref{Meln_nonsym_B}) to eliminate $B$, the identity $\sigma = \lambda (\beta+2)-1$ to eliminate $\sigma$, and the identity (\ref{Meln_nonsym_3}) to eliminate $\lambda$ we obtain 
		\begin{equation}\label{tr_unstable}
			\displaystyle \tr DM = -\frac{4 K k_2 (1+\beta + \sigma)}{\sqrt{B}}<0.
		\end{equation}
		Since $\tr DM < 0$ and $\det DM < 0$, by (\ref{characteristic}), for any sufficiently small $\eps>0$, there is one eigenvalue inside the unit circle and one outside, thus these periodic orbits are of saddle type.
	\end{proof}

	\subsection{Homoclinic orbits ($D_3$)}\label{sec:hom}
	
	In $D_3$, the homoclinic solutions from (\ref{TanhSech}) are given for every $B > 0$ in polar coordinates as 
	\[ \begin{array}{l}
		\displaystyle \xi =  2 \sqrt{B} \sech( \sqrt{B} \tau ) \\
		\displaystyle \sin \phi = - 2 \tanh(\sqrt{B} \tau) \sech(\sqrt{B} \tau)  \\
		\displaystyle \cos \phi = \big ( 1 - 2 \tanh^2(\sqrt{B} \tau) \big ).
	\end{array}
	\]
	For $\eps>0$ the slow evolution of the quantities $A$ and $B$ can be written as (cf.\ \cite{Sparrow}):
	\begin{equation}\label{ABeps}
		\begin{aligned}
			\dot A &= \varepsilon(- \sigma \xi^2 + \beta \zeta + \beta \sigma) = \varepsilon(- \sigma \xi_0^2 + \beta \zeta_0 + \beta \sigma) + O(\varepsilon^2)\\
			B \dot B &= -\varepsilon(\eta^2 + \beta \zeta^2 + \beta \sigma \zeta) =  -\varepsilon(\eta_0^2 + \beta \zeta_0^2 + \beta \sigma \zeta_0) + O(\varepsilon^2),
		\end{aligned}
	\end{equation}
	where zero index denotes the unperturbed solutions from \eqref{TanhSech}.

	As already mentioned in Remark~\ref{rem:hom}, in the present case the only equilibrium point which survives the perturbation  $0<\varepsilon \ll 1$ is $(\xi, \eta, \zeta) = (0, 0, -\sigma)=: \Xi_\sigma$, corresponding to the origin in the original Lorenz system (\ref{Lorenz63}), and we have $B = A = \sigma$ at that point. 
	For $\varepsilon > 0$ the equilibrium $\Xi_\sigma$ has two stable eigenvalues 
	$\displaystyle \nu_1 = -\sqrt{\sigma} + O(\eps)$ 
	and $\nu_2 = - \varepsilon \beta$ and the unstable one 
	$\displaystyle \nu_3 = \sqrt{\sigma} + O(\eps)$. 
	The leading stable direction corresponding to the eigenvalue $\nu_2$ is the invariant line $\{\xi = 0, \; \eta = 0 \}$ of (\ref{Lorenz_resc}), which is $\{A = B\}$ in (\ref{ABeps}).
	
	We note that the second equation of system (\ref{ABeps}) determines the dynamics of variable $B$,  and is in particular equivalent to the last equation of system (\ref{Ham_2}).  Integration over the unperturbed homoclinic orbit gives, to leading order, 
	\begin{equation}\label{e:Bdiff}
		B^2(+\infty) - B^2(-\infty) = -4/3(\beta + 2)\sigma^{3/2}\eps,
	\end{equation}
	which is nonzero for $\eps>0$ since $\beta>0$. 
	This is inconsistent with the persistence result for perturbed homoclinic orbits in \cite{WHhom}, where it is claimed that this integral always vanishes. Hence, \cite{WHhom} is not applicable in the present situation (and is not valid in the claimed generality). 
	
	Towards a correct prediction, first note that a homoclinic orbit converges to the equilibrium point as $\tau \to -\infty$ along the unstable direction, and $\lim\limits_{\tau \to -\infty} A(\tau) = \lim\limits_{\tau \to -\infty} B(\tau) = \sigma$.  Let us assume that a generic homoclinic orbit exists for $0<\eps\ll 1$, which approaches the equilibrium along the leading direction. In the present case this is the line $\{A = B\}$, along which it slowly converges to the equilibrium with rate $\nu_2 = -\varepsilon \beta$. For $0<\eps\ll 1$ this slow part dominates the Melnikov integral, which means a connection of the unstable manifold along $\{A=B\}$ requires that the jumps of $B^2$ in \eqref{e:Bdiff}, and of $A^2$ coincide to leading order in $\eps$. The latter can be computed as
	\begin{equation}\label{ABeps_sol}
		\begin{array}{l}
			\int\limits_{-\infty}^{+\infty}A(\tau)\dot A(\tau)d \tau = 4\varepsilon(\beta - 2 \sigma)\sigma^{3/2}, 
		\end{array}
	\end{equation}
	which equals \eqref{e:Bdiff} if and only if $\sigma = \frac{1}{3}(1 + 2 \beta)$, that is, $\lambda = 2/3$. This is consistent with the discussion of periodic orbits in Proposition \ref{p:per} and, although we do not give a full proof, we thus expect there exists a homoclinic bifurcation curve that tends to $\lambda = 2/3$ when $\varepsilon \to 0$. 
	
	\section{Implications for infinite time averages}
	
	In this section we use the angle bracket notation to denote the infinite time average:
	\[ \langle f(\textbf{X}) \rangle := \lim_{t\to \infty} \frac{1}{t} \int_0^t f(\textbf{X}(s))ds \]
	In \cite{Goluskin2018}, Goluskin illustrates an application of semi-definite programming to dynamical systems by obtaining bounds on time averages for the Lorenz equations.  He proves that the time averages of the following monomials are maximized by the value obtained at the fixed points $\textbf{X}_{\pm}$ for a wide range of parameters $(\beta, \sigma)$ and for all $0 < \rho < \infty$:
	\[ \langle Z \rangle \hspace{.5 cm} \text{ , } \hspace{.5 cm} \langle X^2 \rangle \hspace{.5 cm} \text{ , } \hspace{.5 cm} \langle XY \rangle \hspace{.5 cm} \text{ , } \hspace{.5 cm} \langle Z^2 \rangle \hspace{.5 cm} \text{ , } \hspace{.5 cm} \langle XYZ \rangle \hspace{.5 cm} \text{ , } \hspace{.5 cm} \langle Z^3 \rangle \hspace{.5 cm} \text{ , } \hspace{.5 cm} \langle XYZ^2 \rangle . \]
	Due to the form of the Lorenz equations, certain infinite time averages must be proportional.  For instance, one has
	\[ \langle XY \rangle = \lim_{t\to \infty} \frac{1}{t} \int_0^t X(s)Y(s) ds = \lim_{t\to \infty} \frac{1}{t} \int_0^t \big [ \beta Z(s) + Z'(s) \big ] ds = \beta \langle Z \rangle + \lim_{t\to \infty} \frac{Z(t)-Z(0)}{t} = \beta \langle Z \rangle . \]
	In this way one has the following equalities
	\[ \langle X^2 \rangle =  \langle XY \rangle = \beta \langle Z \rangle \hspace{.5 cm} \text{ , } \hspace{.5 cm} \langle XYZ \rangle = \beta \langle Z^2 \rangle \hspace{.5 cm} \text{ , } \hspace{.5 cm} \langle XYZ^2 \rangle = \beta \langle Z^3 \rangle , \]
	and hence it suffices to consider the following three time averages:
	\[ \langle Z \rangle \hspace{.5 cm} \text{ , } \hspace{.5 cm} \langle Z^2 \rangle \hspace{.5 cm} \text{ , } \hspace{.5 cm} \langle Z^3 \rangle . \]
	
	As stated previously, for $\lambda > 1$ the fixed points become unstable for sufficiently large $\rho$, hence while these sharp upper bounds are indeed valid, they do not represent the values that most trajectories obtain.  However, with the results of the previous section in hand we can provide complementary results which give values for these time averages which are observed for all trajectories within the nontrivial basin of attraction of the symmetric periodic orbit.  For such initial conditions the infinite time averages above are given by the average over the periodic orbit, i.e.,
	\[ \langle Z \rangle =  \lim_{t\to \infty} \frac{1}{t} \int_0^t Z(s,\varepsilon) ds = \frac{1}{ \varepsilon^2 T_1^{\varepsilon}} \int_0^{T_1^{\varepsilon}} \big [ \sigma^{-1} \zeta_1 (\tau',\varepsilon) + 1 \big ] d\tau' , \]
	where $T_1^{\varepsilon}$ is the period of the solution $\zeta_1(\tau, \varepsilon)$.  It seems likely these values do not have a simple closed form and instead we compute the time averages via an expansion in $\varepsilon$.  For instance, the lowest order term can be found from the explicit formulas for the unperturbed solution and period, $\zeta_1(\tau)$ and $T_1$, as follows
	\begin{align*}
		\langle Z \rangle & = \frac{1}{\varepsilon^2} \Big [ 1 + \frac{1}{\sigma T_{1}} \int_0^{T_{1}} \zeta_1 (\tau') d\tau' + \mathscr{O}(\varepsilon) \Big ] = \frac{1}{\varepsilon^2} \Big [ 1 + \frac{\sqrt{B}}{\sigma 4K} \int_0^{\frac{4 K}{\sqrt{B}}} B \big ( 1-2k_1^2 \mathsf{sn}(\sqrt{B}\tau') \big ) d\tau' + \mathscr{O}(\varepsilon) \Big ] \\ &  = \frac{1}{\eps^2}\left [ 1 - \frac{B}{\sigma}\left(1 - \frac{2 E}{K}\right)  +O(\varepsilon) \right ].
	\end{align*}
	In this way, we obtain the following expressions for the infinite time averages, expressed in terms of $\rho$ rather than $\varepsilon$:
	\begin{equation}\label{EpsilonExpansion}
		\begin{split}
			\langle Z \rangle & = \rho \left [ 1 - \frac{B}{\sigma}\left(1 - \frac{2 E}{K} \right) \right ] + \mathscr{O}(\rho^{1/2} ) \\
			\langle Z^2 \rangle & = \rho^2 \Big [ 1 - \frac{3\beta (K-2E)^2}{d_1^2 K} \big ( 8e_1 + \beta e_2 \big ) \Big ] + \mathscr{O}(\rho^{3/2}) \\
			\langle Z^3 \rangle & = \rho^3 \Big [ 1 - \frac{9\beta (K-2E)^2}{5 d_1^3 K} \Big ( 20d_1e_1 + \beta^2 (K-2E) \big ( (32 k^4-36k^2 +19) K - (64 k^4 -64k^2+34 ) E \big ) \Big ) \Big ] + \mathscr{O}(\rho^{5/2}) .
		\end{split}
	\end{equation}
	Recall the expressions $e_1,e_2,d_1$ defined in \eqref{EllipticIntegralExp1}, \eqref{Meln_sym_B} were shown to be positive.  Hence for the time averages $\langle Z \rangle$ and $\langle Z^2 \rangle$, these expressions resolve the coefficient of the leading order term as a function of $\sigma, \beta$ which is strictly less than one, hence less than that of the fixed point value.  For $\langle Z^3 \rangle$, the term $d_1 e_1$ is always positive, whereas one can check that the other expression inside the parentheses is positive for all $\lambda > 2.5611...$, whereas is is negative for $\lambda$ less than this.  Hence by fixing such $\lambda$ and choosing $\beta$ sufficiently large this exceeds the fixed point value.  This agrees with Goluskin's result however, since the region in parameter space where $\langle Z^3 \rangle$ is maximized at $\textbf{X}_{\pm}$ does not include large $\beta$. 
	
	As mentioned previously, the transport $\langle XY \rangle$ is of particular interest, since this is the truncated version of the Nusselt number from fluid dynamics and hence the most well studied such average.  Towards understanding the organization of solution branches and the associated transport, we denote the transport of the stable symmetric periodic orbit by
	\[ H_1 = \beta \langle Z \rangle = \beta \rho \left [ 1 - \frac{B}{\sigma}\left(1 - \frac{2 E}{K} \right) \right ] + \mathscr{O}(\rho^{1/2} ),  \]
	and we also compute the transport obtained by the unstable, asymmetric periodic orbits
	\[	H_2 = \frac{\beta}{\eps^2}\left(1 + \frac{B}{\sigma T_2} \int\limits_0^{ T_2}
	1 - 2  \sn^2 \left( \frac{\sqrt{B}}{k_2} \tau' \right) d\tau' + O(\varepsilon)  \right) 
	= \frac{\beta}{\eps^2}\left(1 - \frac{B}{\sigma}\frac{K(2 - k_2^2) - 2 E}{K k_2^2} + O(\varepsilon)  \right).
	\]
	Next, we cast $H_1, H_2$ in terms of $\rho$ and rescale to a finite range of transport values. This gives
	\begin{align}
		h_{\rho,1}(\lambda)&:= H_1/\rho =  \beta \left(1 -  R_1(\lambda)  + O(\rho^{-1/2}) \right), \; 
		R_1(\lambda) := \frac{B}{\sigma}\left(1 - \frac{2 E}{K}\right), \; \lambda\in(2/3,\infty), \label{e:transport_sym}\\
		h_{\rho,2}(\lambda)&:= H_2/\rho = \beta \left(1 -  R_2(\lambda)  + O(\rho^{-1/2}) \right), \; 
		R_2(\lambda) := \frac{B}{\sigma} \frac{K(2 - k_2^2) - 2 E}{K k_2^2}, \; \lambda\in(2/3,1).\label{e:transport_asym}
	\end{align}
	We first show that $h_{\rho,1}(2/3)=h_{\rho,2}(2/3)=0$ understood as the limit $\lambda \searrow 2/3$, and analogously 
	$h_{\rho,1}(1)= h_{\rho,2}(\infty)=\beta$. Indeed, $R_1(2/3)=R_2(2/3)=1$, $R_1(\infty)=R_2(1)=0$ due to the following.
	\begin{itemize}
		\item $R_1(2/3)=1$: $\lambda\to 2/3$ gives $k_1, k_2\to 1$ so that $E(1)=0$, $K(1)=\infty$ and \eqref{Meln_sym_B}, \eqref{Meln_nonsym_B} imply $B/\sigma\to 1$;
		\item  $R_2(2/3)=1$:  the previous also implies $(K(2 - k_2^2) - 2 E)/(K k_2^2)\to 1$;
		\item $R_1(\infty)=0$: $\lambda\to\infty$ means $k_1\to k_*$, $B$ is bounded and $E/K\to 1/2$;
		\item $R_2(1)=0$: $\lambda\to 1$ gives $k_2\to 0$ and using \eqref{Meln_nonsym_B} as well as $E=K= \pi/2$ at $k_2=0$ implies  
		\[
		R_2(\lambda) = \frac{B}{\sigma} \left(\frac{2(K-E) - K k_2^2}{K k_2^2}\right)
		= \sigma\frac{2(K-E) - K k_2^2}{4 \sigma E + \beta ( 2(K-E) -k_2^2 K)} \to 0 \text{ as $k_2\to 0$.}
		\]
	\end{itemize}
	The differences to the scaled maximum transport $H_{\pm}/\rho= \beta(1-\rho^{-1})$ are the positive quantities
	\begin{equation}
		\beta(1-\rho^{-1})-h_{\rho,j} =\beta(R_j(\lambda)+O(\rho^{-1/2})), \; j=1,2.
	\end{equation}
	In particular, the stable symmetric periodic orbits $L$ yield the same order of magnitude of transport with respect to $\rho$, but feature a $\lambda$ dependent downshift that vanishes at $\lambda=\infty$, i.e. at $\Xb_{\pm}$.  In addition, from our perhaps rough error estimates we obtain a correction of order (at least) $\rho^{-1/2}$ compared to the next order being $\rho^{-1}$ for $\Xb_{\pm}$. Numerical computations suggest that this term might in fact be of order $\rho^{-1}$, but we do not explore this further here.

	\section{Numerical computations and hysteresis loop}\label{sec:numerics}
	
	\begin{figure}
		\begin{center}
			\begin{tabular}{cc}
				\includegraphics[width=0.42\textwidth]{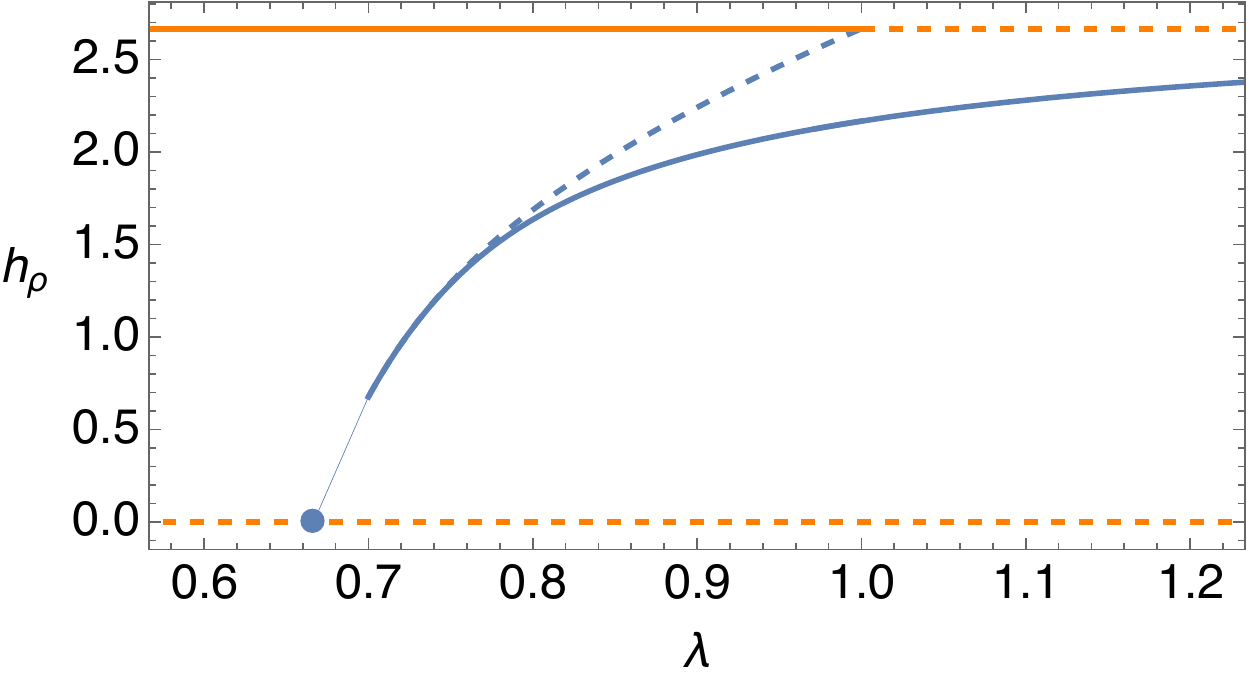}
				&\includegraphics[width=0.42\textwidth]{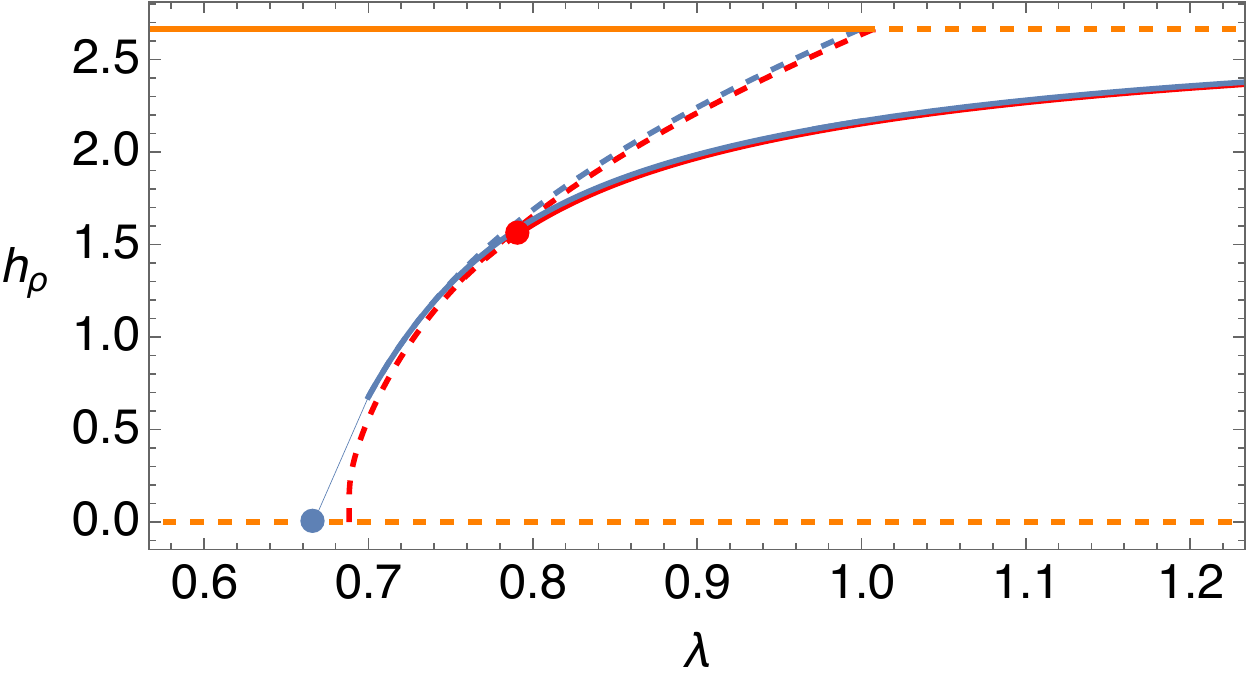}\\
				(a) & (b)
			\end{tabular}
			\caption{Bifurcation diagrams (solid=stable, dashed=unstable) of relevant equilibria and periodic orbits in terms of $\lambda$ and the rescaled transport, illustrating the hysteresis. 
				(a) equilibria in the averaged planar system at $\rho=\infty$ that persist for large finite $\rho$ as equilibria (orange) or periodic orbits (blue). Computations are done via the elliptic integrals with \textsc{Mathematica}. The bullet marks $\lambda=2/3$ at zero transport; for illustration, the thin blue line extends numerical computations to the theoretical limit at zero transport. 
				(b) Overlay of (a) with branches of stable or unstable symmetric (red solid/dashed) and unstable asymmetric (red dashed) periodic solutions to the full system at $\rho=1000$ computed with \textsc{Auto} \cite{auto}. }
			\label{f:hysteresis}
		\end{center}
	\end{figure}
	
\begin{figure}
		\begin{center}
			\begin{tabular}{ccc}
				\includegraphics[height=45mm]{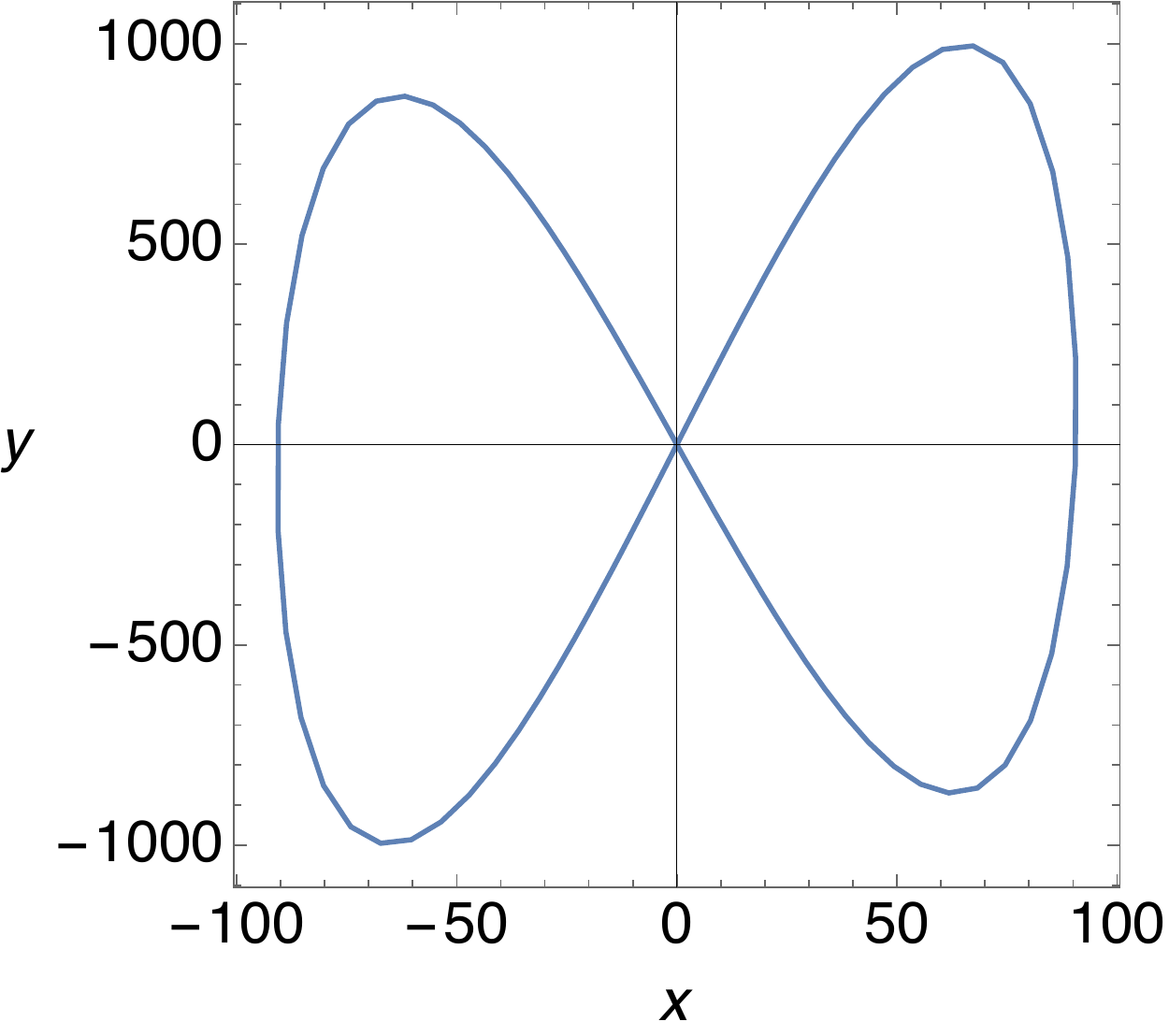}
				&\includegraphics[height=45mm]{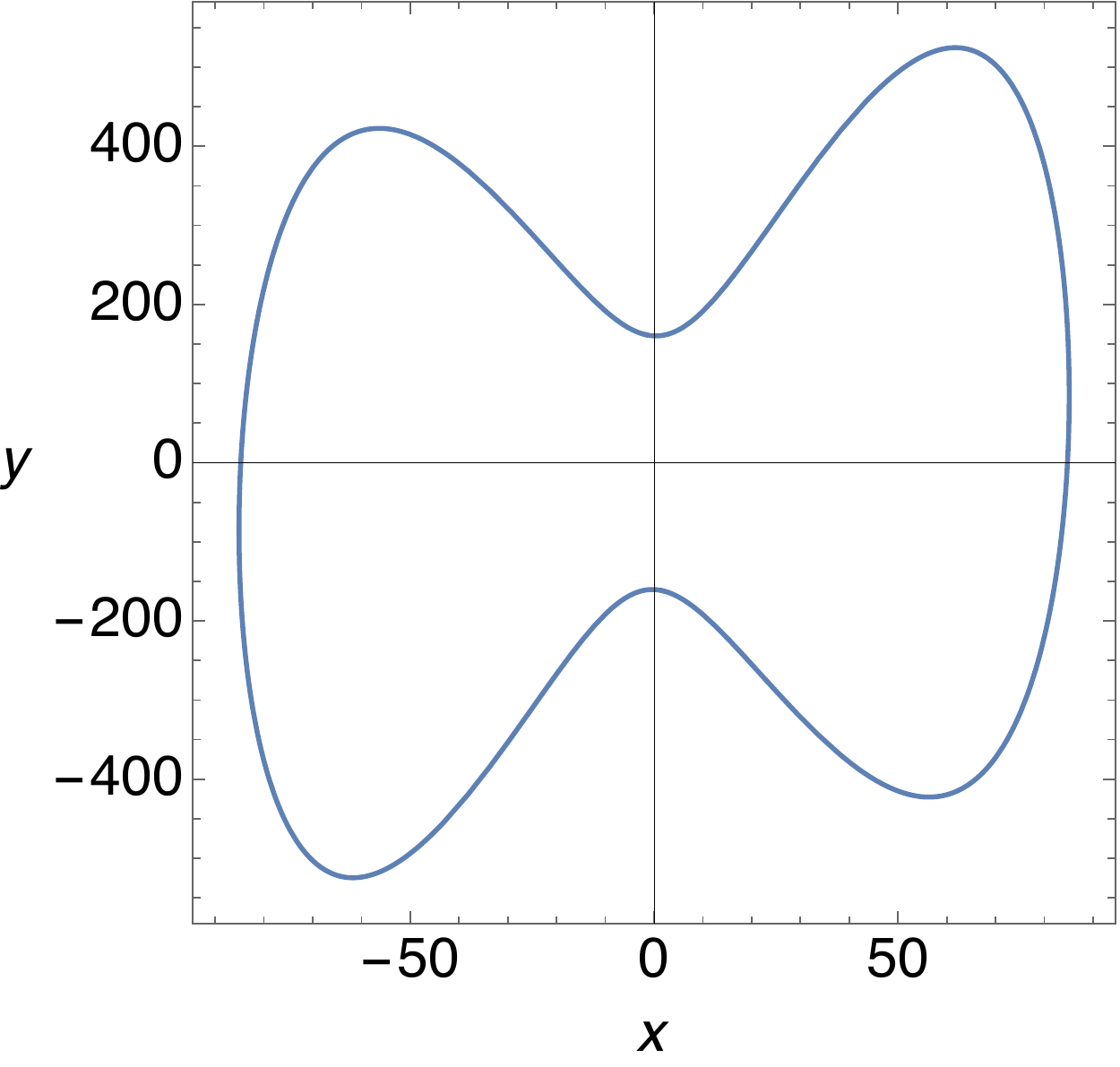}
				&\includegraphics[height=45mm]{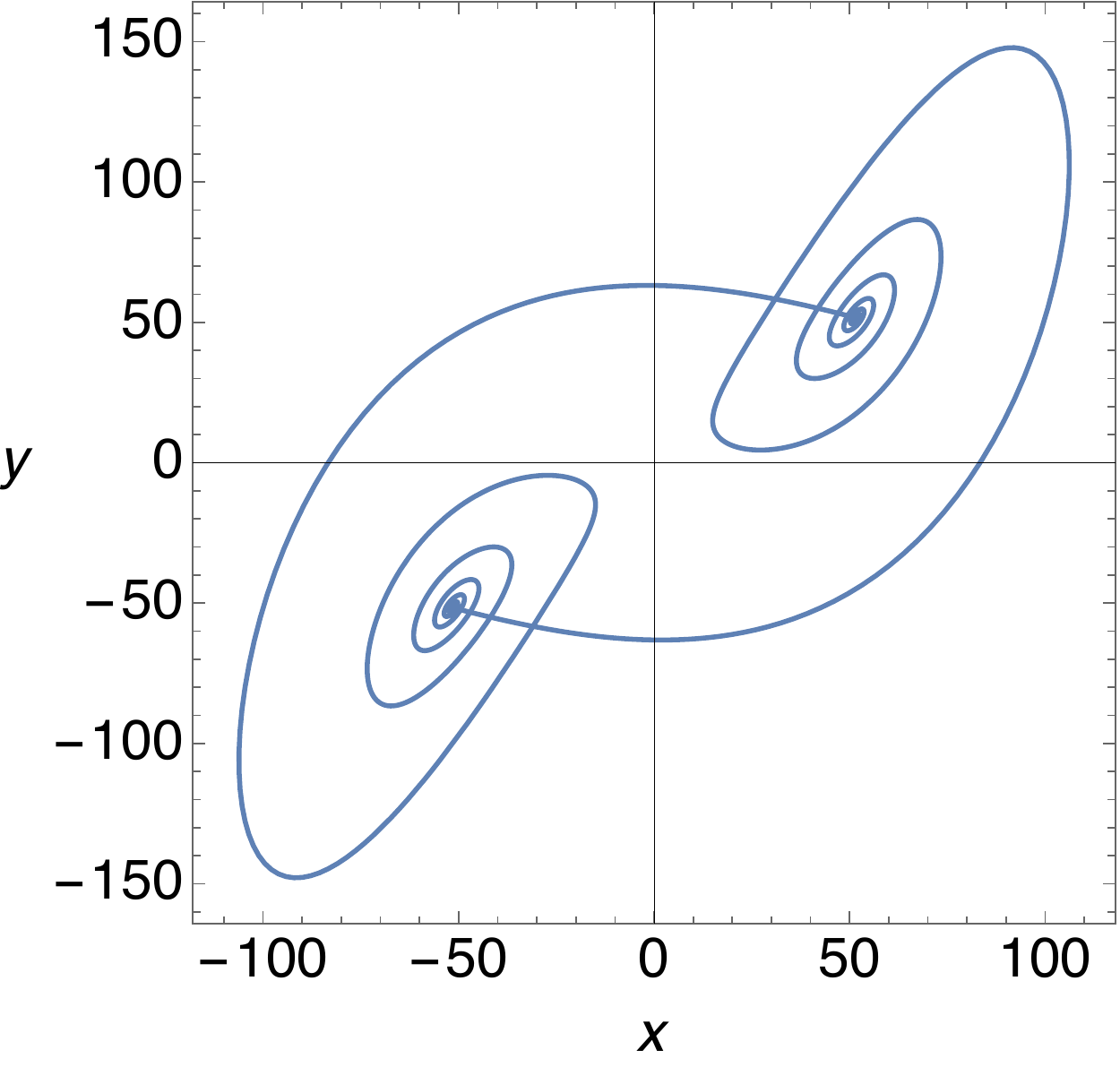}\\
				(a)  & (b) & (c)
			\end{tabular}
			\caption{Profiles of periodic solutions for $\rho=1000$ from the red solid branch in Figure~\ref{f:hysteresis}(b). (a) near the double homoclinic loop at the left termination point; (b) at $\lambda\approx1$; (c) near an apparent symmetric heteroclinic cycle $\lambda\approx 11.3$.}
			\label{f:nearDblHom}
		\end{center}
	\end{figure}

	We present numerical results that corroborate the analytical results for $\rho=\infty$ and $1\ll\rho<\infty$ of the previous section, and that highlight the occurrence of a hysteresis loop. 
	In Figure~\ref{f:hysteresis}(a) we plot the numerical evaluation of \eqref{e:transport_sym} for symmetric periodic orbits and \eqref{e:transport_asym} for asymmetric ones. However, the elliptic integral routines of the current version of the software \textsc{Mathematica} for $B,E,K$ have failed to numerically converge for transport below $\approx 0.6$. The analytical prediction is that the branches terminate at $\lambda=2/3$ in homoclinic bifurcations of the zero equilibrium and thus at zero transport. Indeed, at $\lambda=2/3$ the intersection of the level sets of the conserved quantities $A, B$ forms a symmetric pair of homoclinic loops, cf.\ Fig.~\ref{f:nearDblHom}(a), which is the limit of the branch of symmetric periodic orbits, and each branch of asymmetric periodic orbits limits on one of the homoclinic loops. 
	
	The arrangement of branches in Figure~\ref{f:hysteresis}(a) together with the stability properties suggest a hysteresis loop of equilibria and periodic orbits in terms of $\lambda$: For $\lambda<2/3$ the equilibria $\textbf{X}_{\pm}$ that maximize transport are stable, while for $\lambda>1$ the  symmetric periodic orbit are. Intermediate values $2/3<\lambda<1$ lie in the analytically predicted region of bistability with stable equilibria $\textbf{X}_{\pm}$ and stable symmetric periodic orbit $L$. 
	
\begin{figure}
		\begin{center}
			\begin{tabular}{cc}
				\includegraphics[height=40mm]{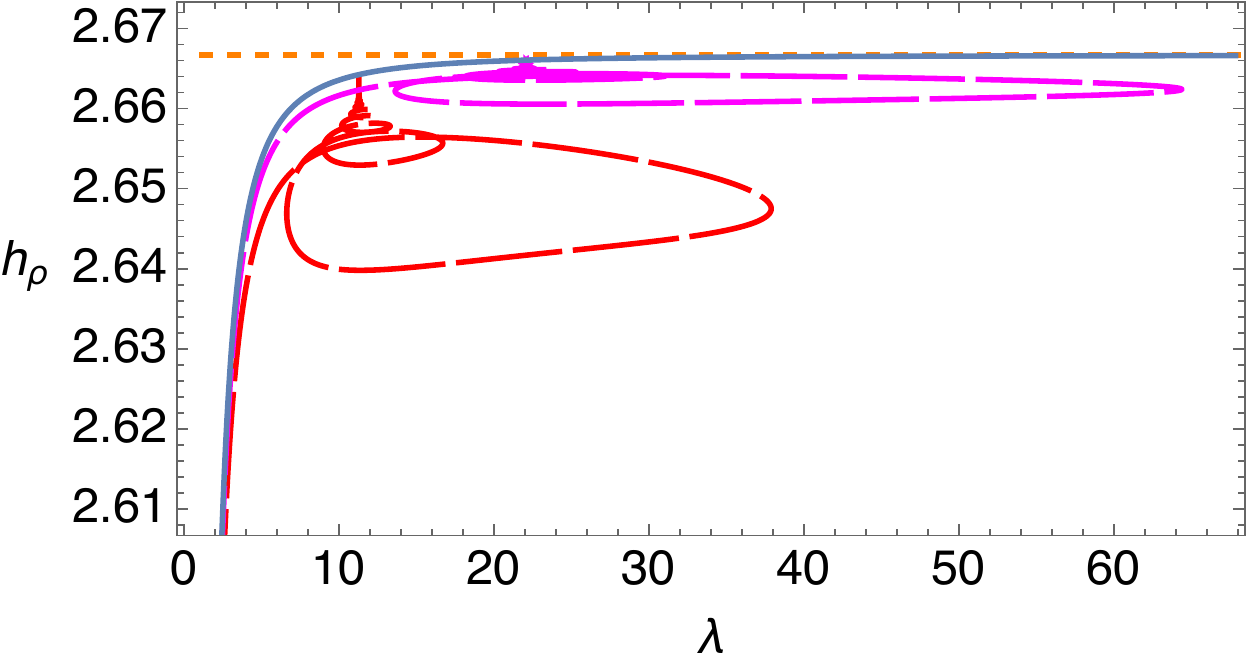}
				&\includegraphics[height=40mm]{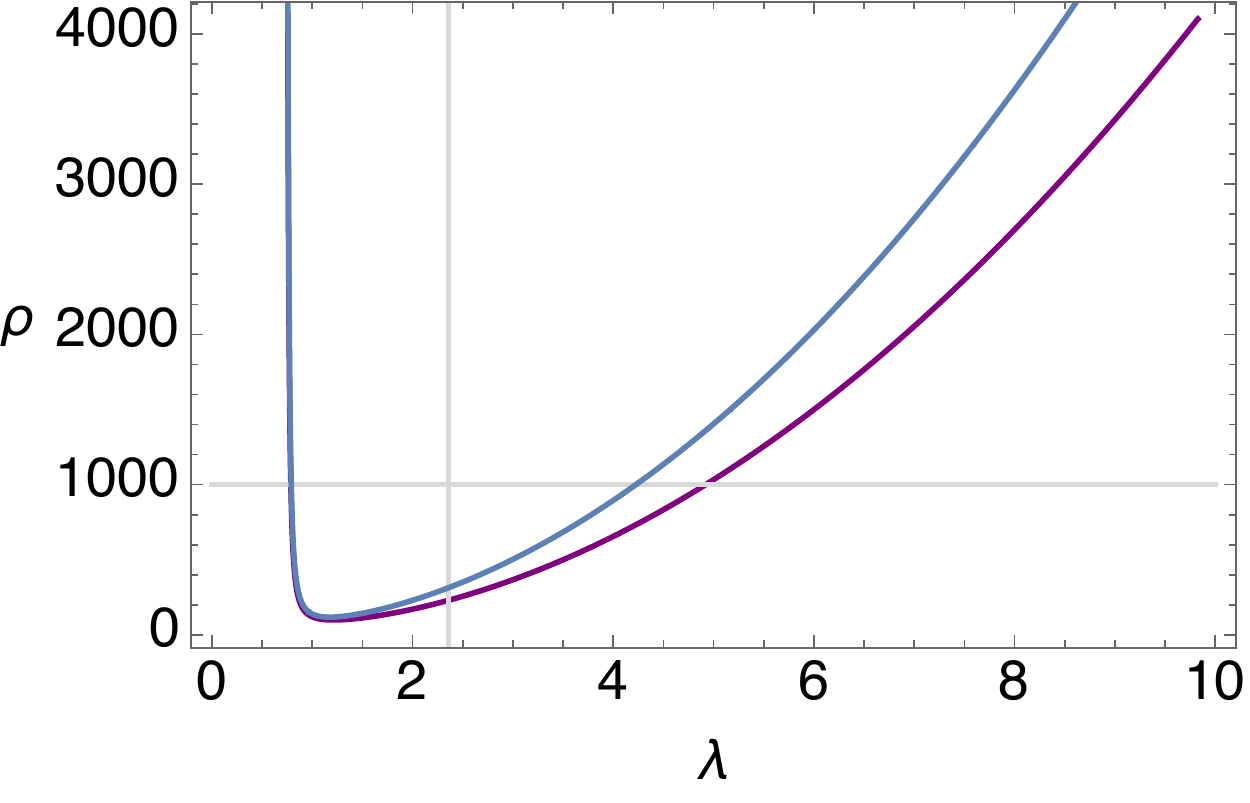}\\ 
				(a) & (b)
			\end{tabular}
			\caption{(a) Bifurcation diagram for larger values of $\lambda$ analogous to right panel of Figure~\ref{f:hysteresis} but without stability information of the periodic solutions (long dashed lines). The red curve corresponds to the extension of the branch of stable symmetric periodic orbits from Figure~\ref{f:hysteresis}. The magenta curve is the analogue for $\rho=4000$. Continuing along these branches from their lower left ends, numerically before each fold a destabilising period doubling bifurcation occurs and the solution restabilises at the fold point. Each branch appears to terminate in a symmetric heteroclinic cycle between the pair of equilibria $\textbf{X}_{\pm}$. (b) loci of branch points of symmetric periodic orbits (blue) and period doubling points on the bifurcating branches (purple); for values of $\rho$ above the blue curve periodic orbits appear to be stable. Gray lines mark $\rho=1000$ (horizontal) and $\lambda\approx 2.36$ at the classical Lorenz values $\sigma=10, \beta=8/3$.}
			\label{f:largelam}
		\end{center}
	\end{figure}

	For large finite $\rho$ and moderate values of $\lambda$, numerical pathfollowing computations using \textsc{Auto} corroborate that branches of symmetric and asymmetric periodic orbits persist as predicted. See Figure~\ref{f:hysteresis}(b). Towards zero transport the branches of symmetric and asymmetric periodic orbits appear to terminate in homoclinic bifurcations to the zero equilibrium near $\lambda=0.688$. See also Fig.~\ref{f:nearDblHom}(a,b). The asymmetric periodic orbits are unstable as predicted, but the symmetric periodic orbits lose stability at low transport. For $\rho=1000$ this occurs at $\lambda=\lambda_{\rm bp}(\rho)\approx 0.79$ in a supercritical pitchfork bifurcation. A branch of stable periodic orbits bifurcates, which are asymmetric in a different sense, but these lose stability at $\lambda=\lambda_{\rm pd}(\rho)\approx 0.787$ in a period doubling bifurcation. We plot the loci of $\lambda_{\rm bp}, \lambda_{\rm pd}$ in Figure~\ref{f:largelam}(b), showing that as $\rho$ increases, $\lambda_{\rm bp}, \lambda_{\rm pd}$ approach $\lambda=2/3$.

	\medskip
	Further numerical simulations corroborate the hysteresis-type loop: for $\lambda<2/3$ the maximum transport equilibria $\textbf{X}_{\pm}$ appear to be global attractors, while for $\lambda>1$ this seems to be the stable symmetric periodic orbit, as in Fig.~\ref{f:nearDblHom}(b). See also Fig.~\ref{f:SouzaDoeringBound}(b). For $2/3<\lambda<1$ the situation with large finite $\rho$ is complicated by the fact that symmetric periodic orbits are born in a homoclinic bifurcation at some $\lambda_{\rm hom}\in(2/3,1)$ and, as mentioned, are unstable until a bifurcation point $\lambda_{\rm bp} \in (\lambda_{\rm hom},1)$. Up to the aforementioned region of stable asymmetric periodic orbits that bifurcate from $\lambda_{\rm bp}$, the global attractors for $\lambda<\lambda_{\rm bp}$ seem to be $\Xb_{\pm}$ and the region of bistability with the symmetric periodic orbit is effectively $\lambda\in(\lambda_{\rm bp}, 1)$. 
	
\begin{figure}
		\begin{center}
			\begin{tabular}{cc}
				\includegraphics[width=0.42\textwidth]{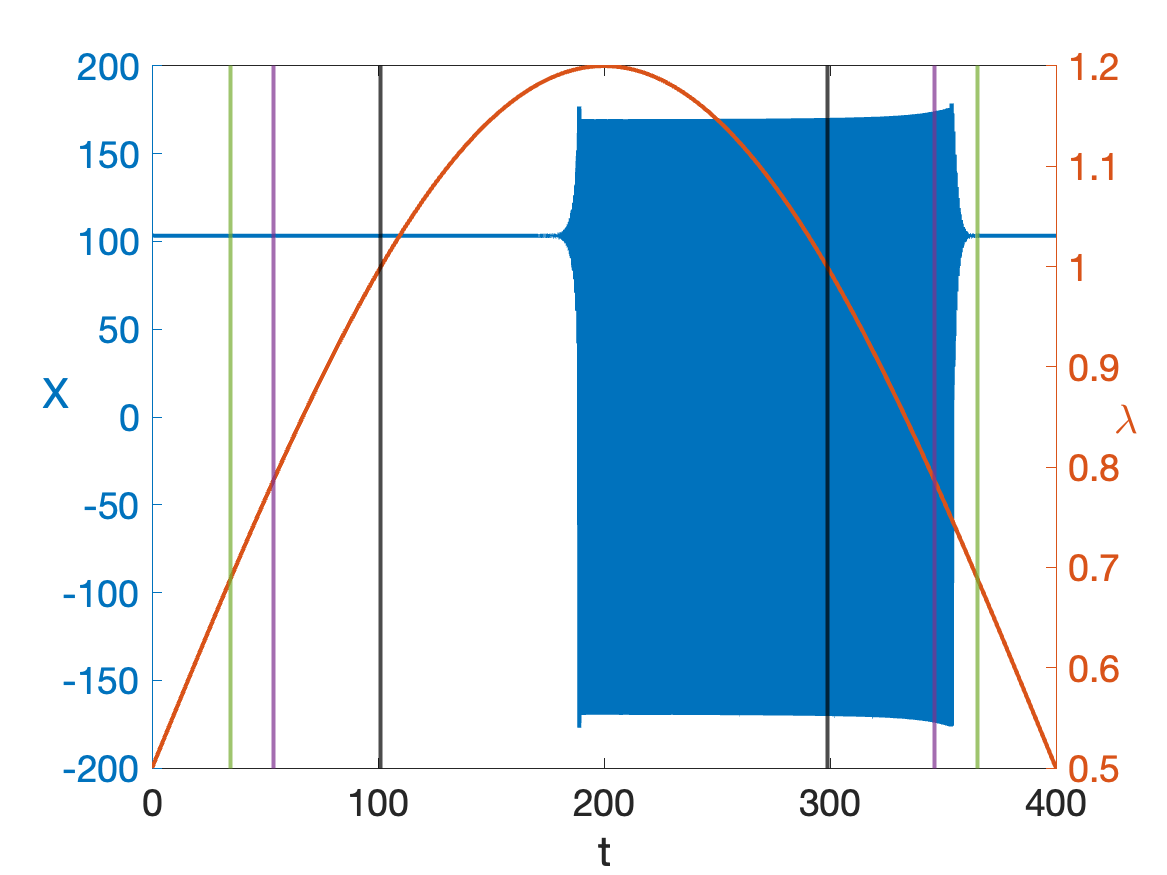}
				& \includegraphics[width=0.42\textwidth]{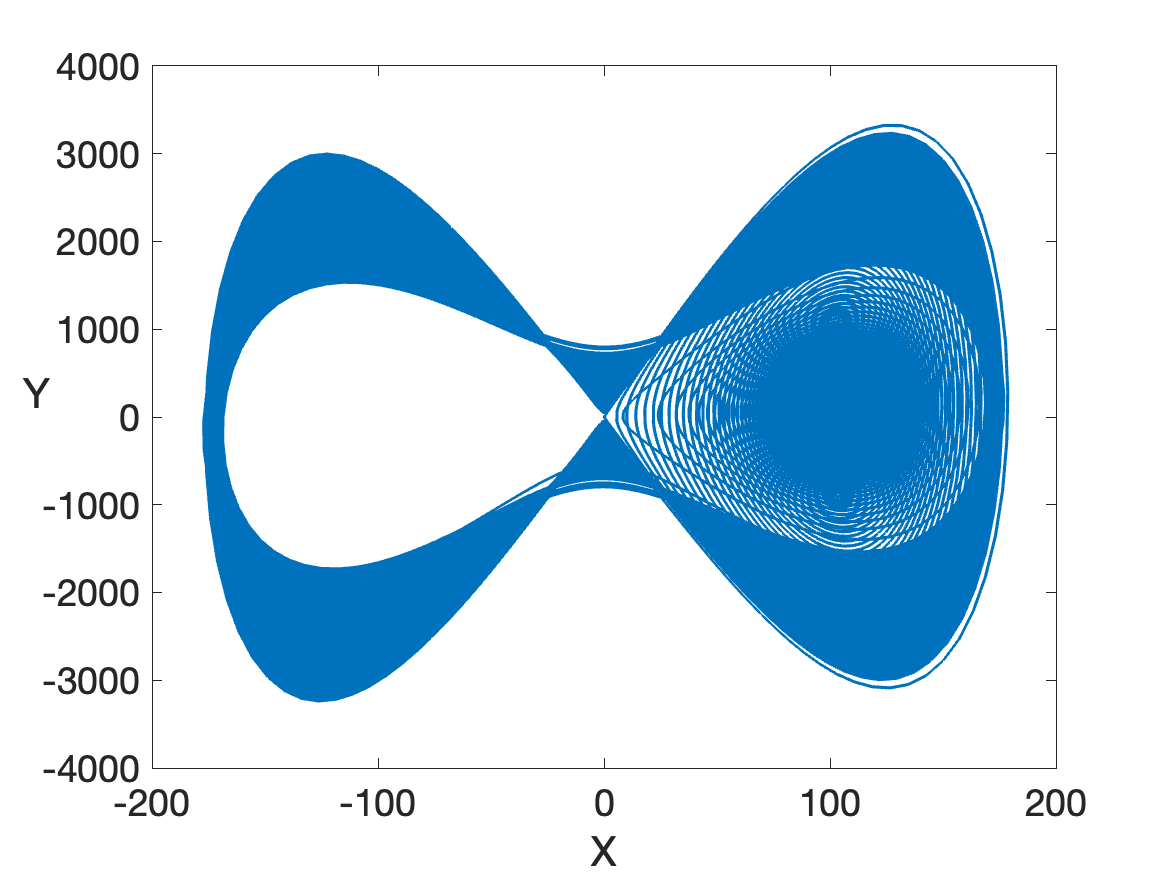}\\
				(a) & (b)
			\end{tabular}
			\caption{We plot a simulation of the hysteresis loop with time-varying $\lambda$ for $\rho=1000$ using Matlab's ode45 routine. (a) The $X$-coordinates of the resulting solution (blue curve, left axis) from a parabolic variation of $\lambda$ (orange curve, right axis). The vertical bars mark the homoclinic bifurcation point $\lambda\approx0.688$ (green), the period doubling $\lambda\approx0.787$ (purple), and the Hopf bifurcation $\lambda=1$ (black). (b) The value of $X$ vs.\ $Y$ of the simulation in (a).}
			\label{f:hysteresissim}
		\end{center}
	\end{figure}

	Using time-varying values of $\lambda$ we consistently found hysteresis as plotted in Figure~\ref{f:hysteresissim}(a): for slowly increasing $\lambda$ from $0$, the solution is quickly close to $\textbf{X}_+$ so that maximum local transport is realized, i.e., transport computed over a time interval of finite length, which can be chosen longer for slower change of $\lambda$. As $\lambda$ increases beyond $1$, the solution eventually approaches the stable symmetric periodic orbit, so that the realized local transport is smaller than the theoretical maximum. Analogous to delayed bifurcations, this transition to the periodic orbit does not occur immediately after crossing $\lambda=1$ at $t=100$, but with a delay, here until around $t=190$. Subsequent decrease of $\lambda$ causes the solution to track the stable branch of symmetric periodic orbits, cf.\ Figure~\ref{f:hysteresissim}(b), which decreases the observed local   transport further until $\lambda=\lambda_{\rm bp} \approx \lambda_{\rm pd}$. Upon decreasing $\lambda$ below this threshold, a switch to a stable equilibrium $\Xb_{\pm}$ occurs, thus re-creating maximum local transport.
	
	\bigskip
	While the asymptotically predicted branch of symmetric periodic orbits of Fig.~\ref{f:hysteresis} continues for increasing $\lambda$ monotonically and unboundedly, we found that for finite $\rho$ this is not the case. As plotted in Figure~\ref{f:largelam}, the branch of stable symmetric periodic orbits turns around, oscillates, and appears to terminate in a symmetric heteroclinic bifurcation of $\textbf{X}_{\pm}$ at a finite value of $\lambda$. See also Fig.~\ref{f:nearDblHom}(c). Upon increasing $\rho$, this turning and termination occurs at larger values of $\lambda$. Hence, this scenario is consistent with the analytical results, which concern $\rho\to \infty$ for bounded ranges of $\lambda$. The appearance of a symmetric heteroclinic cycle between $\textbf{X}_{\pm}$ in the Lorenz system has already been noticed in \cite{Sparrow, GlenSpar}, albeit apparently not in the regime of large $\rho$. 
	
	The transport at such a heteroclinic cycle  is that of the symmetric equilibria, i.e. $\beta(1-\rho^{-1})$, which is indeed very closely matched at the numerical termination points. The $(\lambda,h_\rho)$-loci of the termination points lie near the curve of symmetric periodic orbits (blue solid), which therefore appear to predict the loci of the heteroclinic cycles. The oscillating stability along the branch creates multi-stable regions in $\lambda$; we note that generic unfoldings of the type of heteroclinic cycle with leading oscillating dynamics yield chaotic attractors \cite{bykov_1999}.
	
\begin{figure}[htbp]
		\begin{center}
			\begin{tabular}{cccc}
				\includegraphics[height=40mm]{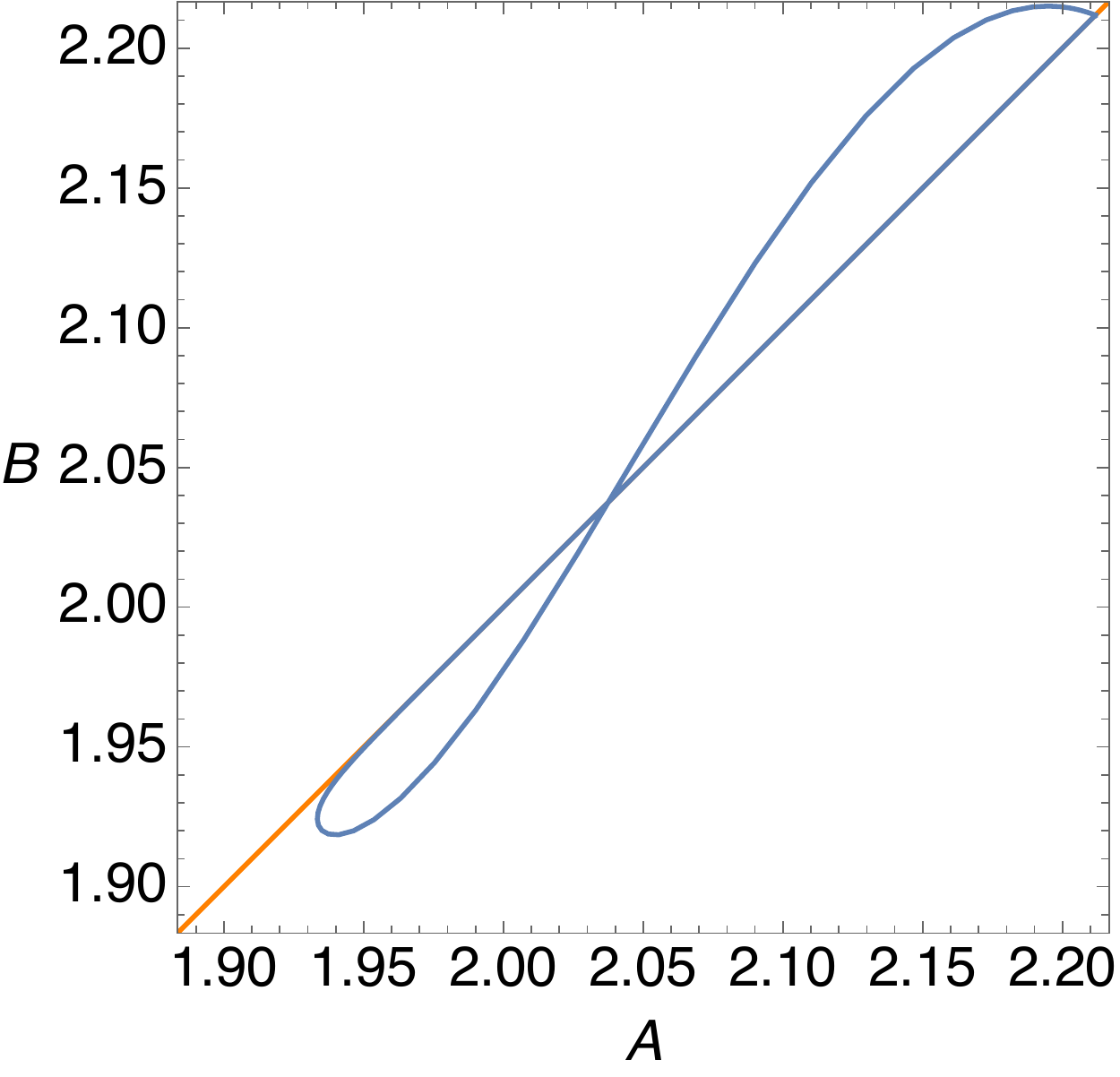}
				&\includegraphics[height=40mm]{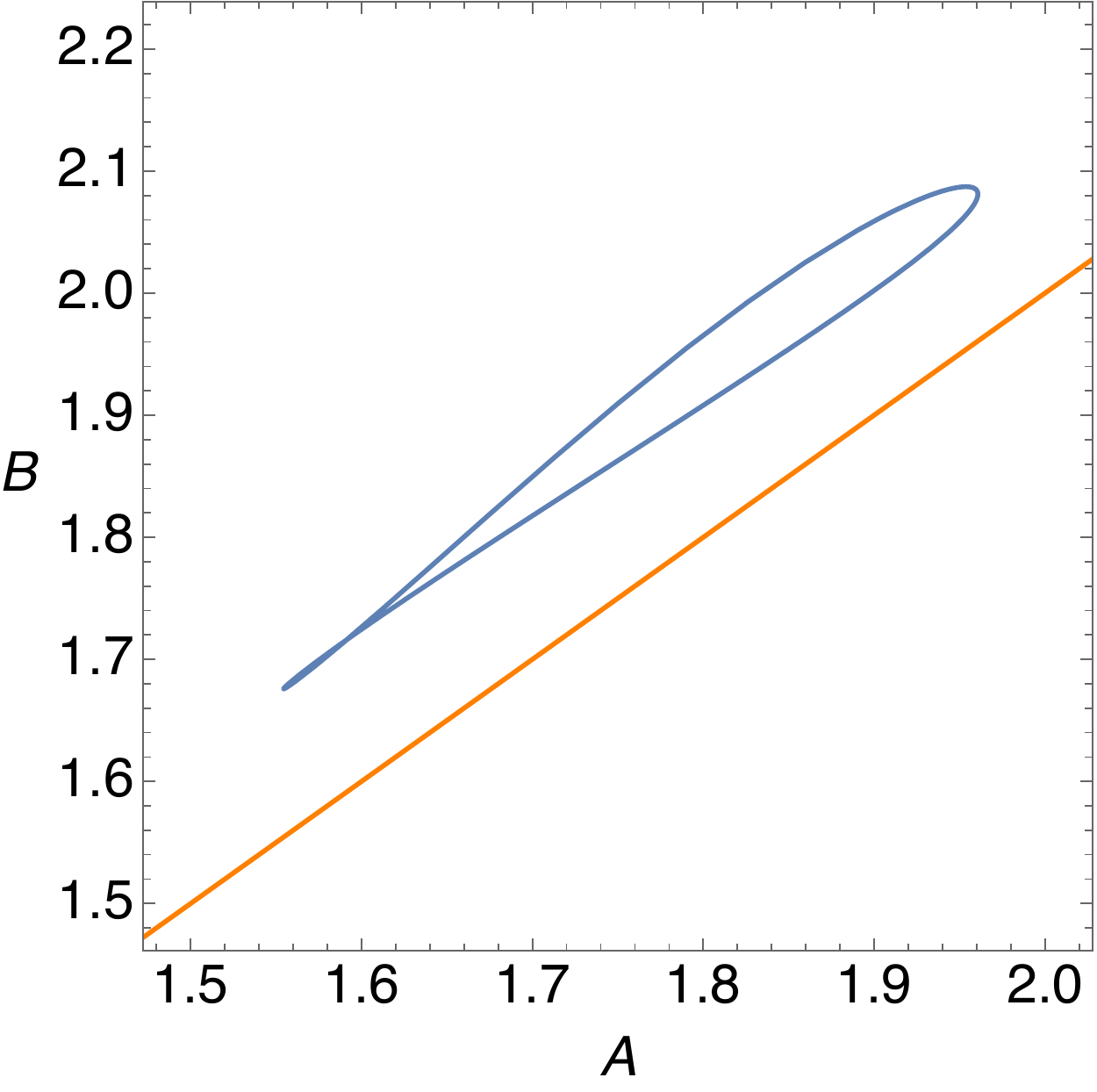}
				&\includegraphics[height=40mm]{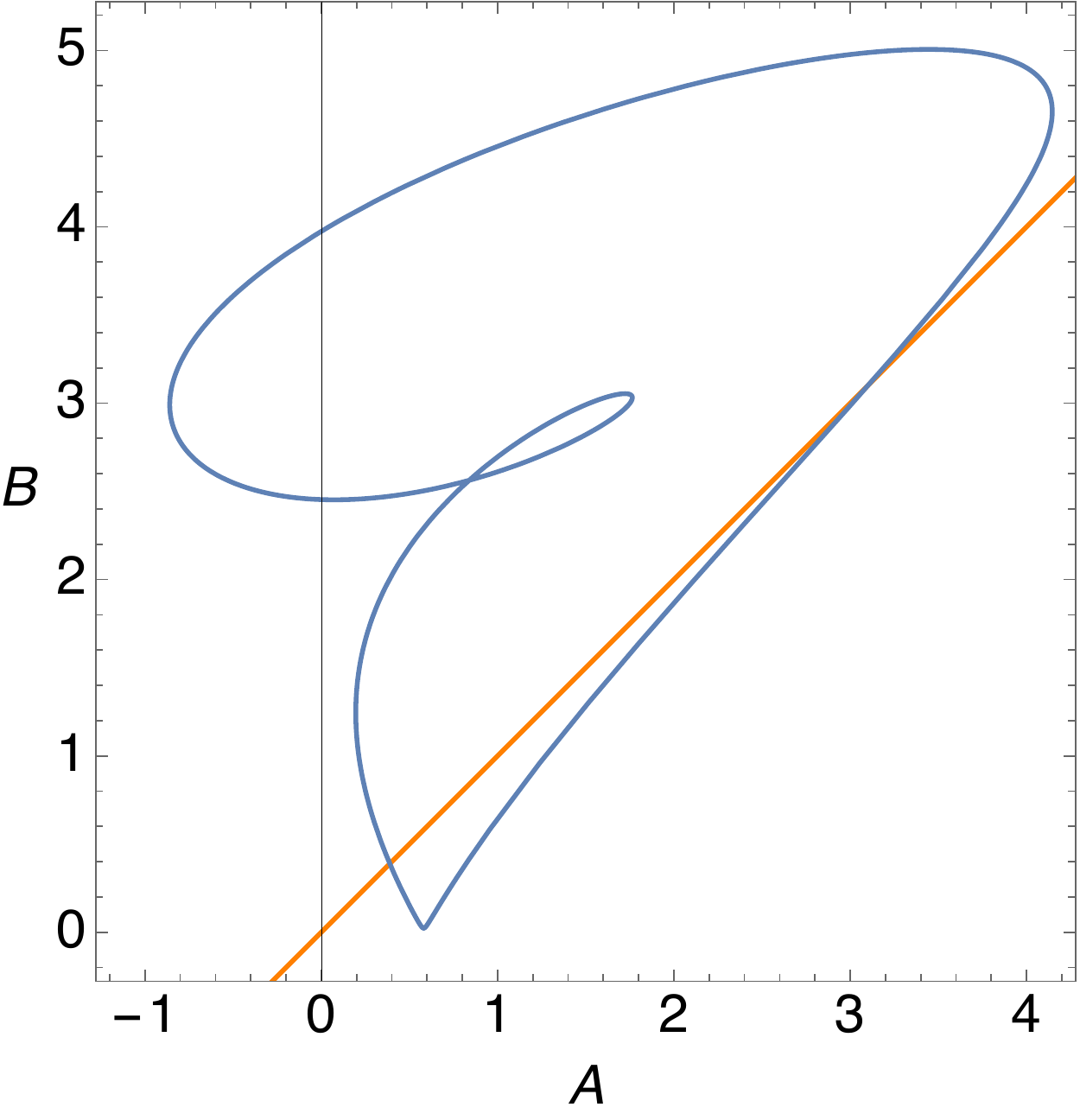}
				&\includegraphics[height=40mm]{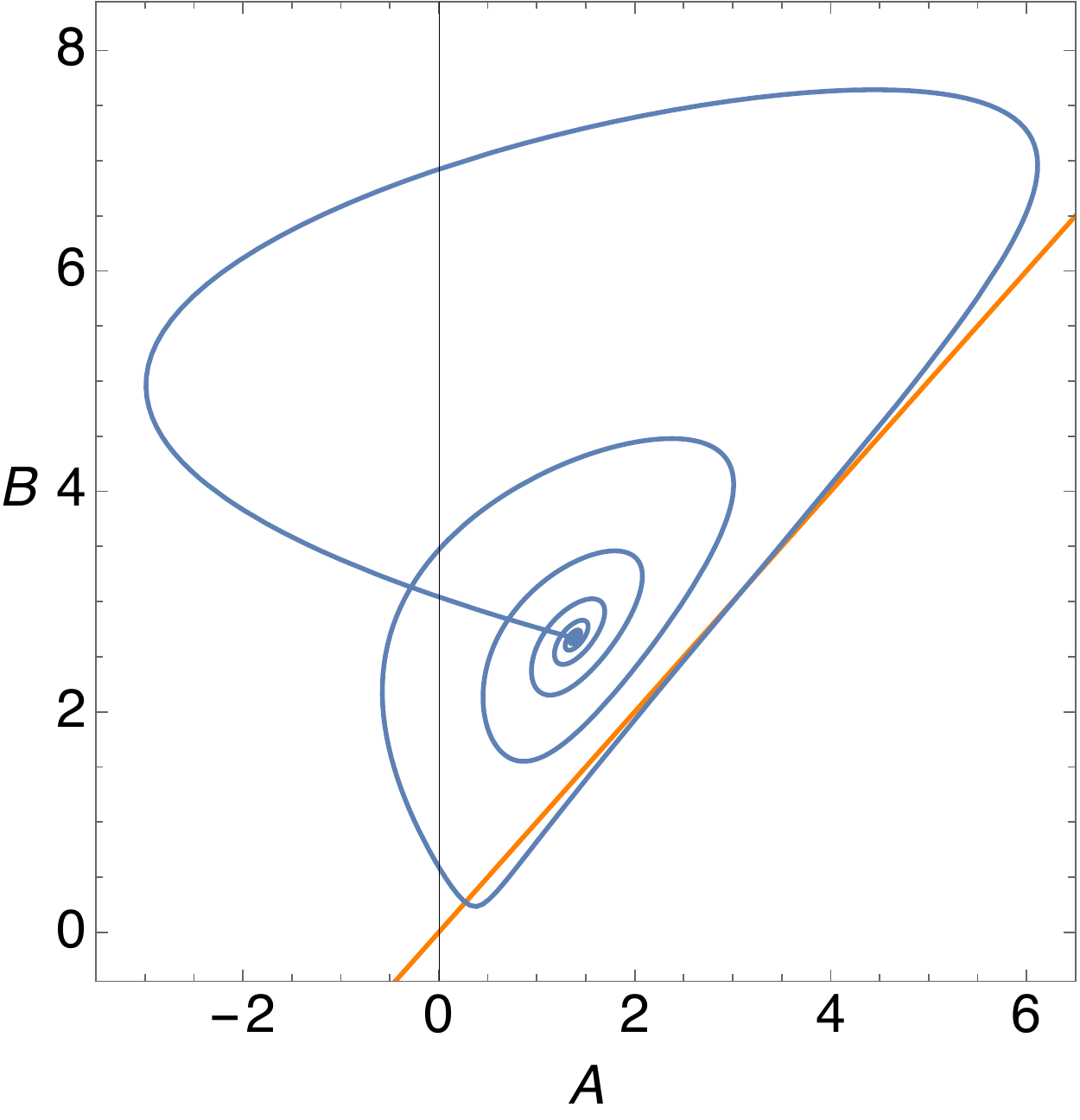}\\ 
				(a) & (b) &(c)&(d)
			\end{tabular}
			\caption{Projections into the $(A,B)$-coordinate plane of periodic solutions for $\rho=1000$ (blue) and the diagonal (orange). (a) Near the double homoclinic loop at the left termination point of the red solid branch in Figure~\ref{f:hysteresis}(b).
				(b) From the red solid branch of Figure~\ref{f:hysteresis}(b) at $\lambda\approx1$. 
				(c) From the red long dashed branch of Figure~\ref{f:largelam}(a) at $\lambda=7.2$, towards the heteroclinic cycle. 
				(d) Near the heteroclinic cycle from Figure~\ref{f:hysteresis}(b).  }
			\label{f:ABproj}
		\end{center}
	\end{figure}

	We plot the projection of a solution near the symmetric heteroclinic cycle into the $(A,B)$-plane in Figure~\ref{f:ABproj}. This corroborates the conjecture by Sparrow in \cite{Sparrow} that orbits bifurcate from $\rho=\infty$ which cross through the diagonal $A=B$. We find that also the solutions near the double homoclinic loop with small transport cross the diagonal. In contrast, the solutions for moderate transport remain in $D_1$ as predicted by the limit $\rho\to\infty$.
	
	\section{Other Lorenz-like systems}

	The analysis for large $\rho$ carries over to other models related to the Lorenz equations \eqref{Lorenz63}. For the general context of extensions, we refer to \cite{Curry1978,Sparrow,Park2021,OlsonDoering2022} and the references therein.  For illustration purposes, let us consider linear additions to \eqref{Lorenz63} in the form
	\begin{equation}
		\begin{aligned}
			\label{LorenzLinAdditons}
			\Xb' &= \Fb(\Xb) + \Ab w + \bb\\ 
			w' &= \Bb (\Xb, w),
		\end{aligned}
	\end{equation}
	with $w\in\mathbb{R}^k$, linear $\Ab, \Bb$ and constant $\bb$. Upon rescaling as in \eqref{rescaling} and $w=\eps^{-j}\omega$, with $\Xib=(\xi,\eta,\zeta)^\intercal$ we obtain the form
	\begin{equation}
		\begin{aligned}
			\dot\Xib &= \Fib_\eps(\Xib) + \mathrm{diag}(\eps^{2-j},\eps^{3-j},\eps^{3-j}) \Ab \omega + \mathrm{diag}(\eps^2,\eps^3,\eps^3)\bb\\ 
			\dot \omega &= \Bb( \mathrm{diag}(\eps^j,\eps^{j-1},\eps^{j-1})\Xib, \eps\omega),
		\end{aligned}
	\end{equation}
	where $\Fib_\eps$ is the right hand side in \eqref{Lorenz_resc}. 
	For $\Ab=\Bb=0$, i.e., in absence of $\omega$, the difference to \eqref{Lorenz_resc} is of order $\eps^2$. Hence,  the leading order analysis is unchanged, which means that periodic orbits bifurcate/persist as for $\bb=0$, although their symmetry properties may be broken. In particular, this applies to the Lorenz models with offsets from \cite{Agarwal,Palmer98} for which one can also show that the transport is maximized in an equilibrium \cite{PreprintIvan}. 
	
	Non-zero $\Bb$ generally requires $j\geq 1$ for a regular limit in which the right hand side of the equation for $\omega$ becomes independent of $\omega$, and vanishes for $j>1$. For $j=1$ we obtain, up to terms of order $\eps^2$, 
	\begin{equation}\label{e:LorenzExtended}
		\begin{aligned}
			\dot\Xib &= \Fib_\eps(\Xib) + (\eps \Ab_1 \omega,0,0)^\intercal \\ 
			\dot \omega &= \Bb( \mathrm{diag}(\eps,1,1)\Xib, \eps\omega),
		\end{aligned}
	\end{equation}
	with $\Ab_1$ the first row of $\Ab$. For $\Bb$ of the form $\Bb = [B_1 | 0 | 0 | B_2]$ the equation for $\omega$ has the slow form  
	\begin{align}\label{e:specialB}
		\dot \omega = \eps (B_1 \xi + B_2\omega), 
	\end{align}
	which occurs with $w\in\mathbb{R}$ in the Lorenz-Stenflo model from \cite{LStenflo_1996}, its magnetic variant \cite{Wawrzaszek}, and with $w\in \mathbb{R}^2$ in the models from \cite{Molteni}; for the latter we  choose $\alpha=O(\eps)$ and shift the auxiliary variables (which gives $\bb\neq 0$) to obtain the form \eqref{e:specialB}. An extension of Lorenz-Stenflo with nonlinear additional equations is considered in \cite{Moon2021}, but still fits into the present framework when, e.g., scaling the variables in addition to Lorenz-Stenflo with $j=3$ and choosing Lewis number of order $\eps^{-2}$. Other extensions of the Lorenz model with two nonlinear auxiliary equations are studied in \cite{costa_knobloch_weiss_1981,Shen2014,Felicio_2018}, which also fit into the present framework when suitably scaling the auxiliary modes and parameters. However, in many cases the situation is more complicated, for instance for the three-dimensional extension in \cite{Shen2015,Felicio_2018}. 
	
	We next show that for the case \eqref{e:specialB} the results of the previous sections also carry over; the following analysis is more explicit in \S\ref{s:Stenflo} for the model from \cite{LStenflo_1996}. In the case \eqref{e:specialB} the Melnikov analysis of \S\ref{s:Melnikov} can be simply extended by adding the slow equation for $\omega$ to the action-angle formulation. The additional Melnikov-integral term $\Mb_3$ is then simply the integral of $B_1 \xi(t) + B_2\omega_0$ over the period $T$, with $\omega_0$ constant. Since $\xi$ has zero average ($\dot\varphi=\xi$ at $\eps=0$ in \eqref{Ham_1}), this term becomes $\omega_0/T$, with the period $T$, so that  $\Mb_3=0$ requires $\omega_0=0$. This means that the values of the other two Melnikov-integrals for \eqref{e:LorenzExtended}, $\tilde M_1, \tilde M_2$, actually coincide with those of $M_1, M_2$ from \S\ref{s:Melnikov}. The non-degeneracy condition turns into invertibility of the matrix  
	\[ DM = \frac{\partial (\tilde M_1,\tilde M_2, \Mb_3)}{\partial (I,B, \omega)}, \]
	where $\tilde M_2$ is independent of $\omega$, and it turns out that also $\tilde M_1$ is: In its integrand $F$ from \eqref{action-angle}, the additional term from $\Ab_1 \omega_0$ is constant and has a factor $\frac{\partial I}{\partial \xi} = \frac{\partial I}{\partial A} \frac{\partial A}{\partial \xi} = T \xi$, where $\xi$ has zero average as noted above. Hence, the matrix has lower left triangular block structure and the block $\partial_\omega \Mb_3 = T B_2$ is invertible, if $B_2$ is. In that case the non-degeneracy condition is therefore the same as for $M_1, M_2$ from the original Lorenz system.

	\subsection{Lorenz-Stenflo}\label{s:Stenflo}
	
	The Lorenz-Stenflo system is given as follows:
	\begin{equation} \label{LorenzStenflo} \begin{split}
			X' & = \sigma ( Y-X ) +s V \\
			Y' & = \rho X - Y - XZ \\
			Z' & = -\beta Z + XY \\
			V' & = - X - \sigma V.
	\end{split}  \end{equation}
	This system is a mode truncation of the rotating Boussinesq equations:
	\begin{equation} \label{RotatingBoussinesqEquations} \begin{split}
			\partial_t \textbf{u} + (\textbf{u} \cdot \nabla ) \textbf{u} + \nabla P  + 2\Omega \hat{z} \times \textbf{u} & = \nu_m \Delta \textbf{u} + \alpha g T \hat{z} \\ 
			\partial_t T + \textbf{u} \cdot \nabla T & = \nu_T \Delta T \\
			\nabla \cdot \textbf{u} & = 0,  \end{split}
	\end{equation}
	where one is considering convection in a fluid in a rotating frame, and a term representing the Coriolis force has been added.  One obtains (\ref{LorenzStenflo}) by making the analogous reduction to a system of ODE's for the Fourier coefficients, but one must include an additional Fourier coefficient $V(t)$ in the expansion of the velocity, which couples to the $X$-mode via the Coriolis force.  The parameter $s$ measures the speed of the rotation.
	
	Since $X$ and $V$ both represent velocity variables, we expect they have the same scaling in $r$, hence we scale
	\[ \epsilon = \rho^{-1/2} \text{ , } X = \epsilon^{-1} \xi \text{ , } Y = \epsilon^{-2} \sigma^{-1} \eta \text{ , } Z = \epsilon^{-2} ( \sigma^{-1} \zeta + 1 ) \text{ , }  V = \epsilon^{-1} \chi \text{ , }   t = \epsilon \tau  \]
	and we obtain the system of equations
	\begin{equation}
		\label{2022_03_PhysExtLorenz_LargeReynRescale}	    \begin{split}
			\frac{d \xi}{d\tau} & = \eta - \epsilon ( \sigma \xi - s \chi ) \\ 
			\frac{d \eta}{d\tau} & = -\xi \zeta - \epsilon \eta \\ 
			\frac{d \zeta}{d\tau} & = \xi \eta - \epsilon \beta (\zeta + \sigma ) \\ 
			\frac{d \chi}{d\tau} & = -\epsilon (  \xi + \sigma \chi ).
		\end{split}
	\end{equation}
	The limiting system when $\epsilon = 0$ coincides with \eqref{Lorenz_resc} except trivial dynamics in the variable $\chi$, so that the system now admits three invariants of motion
	\begin{equation}
		\label{LargeRInvariants}
		\xi^2 - 2\zeta = 2A,\quad \eta^2 + \zeta^2 = B^2,\quad \chi.
	\end{equation}
	Using the first two invariants as for the Lorenz system, (\ref{2022_03_PhysExtLorenz_LargeReynRescale}) can be solved at $\eps=0$, where the solutions have a different form depending on the choice of $(A,B)$ as described in \S \ref{s:Melnikov}.  Analogous to \S\ref{s:Melnikov}, we change coordinates via
	\[ \zeta = B \cos ( \phi ) \hspace{.25 cm} \text{ , } \hspace{.25 cm} \eta = B \sin ( \phi) \hspace{.25 cm} \text{ , } \hspace{.25 cm} \xi = \xi \hspace{.25 cm} \text{ , } \hspace{.25 cm} \chi = \chi \]
	and \eqref{2022_03_PhysExtLorenz_LargeReynRescale} becomes 
	\begin{equation*}
		\begin{array}{lll}
			\dot{\xi} = f_1 + \epsilon g_1 \hspace{2 cm} &  \hspace{4 cm} & g_1 = -  ( \sigma \xi - s \chi )\\
			\dot{\phi} = f_2 + \epsilon g_2 & f_1 = B \sin (\phi)  & g_2 = \sin (\phi) \big [ (\beta-1) \cos (\phi ) + \frac{ \beta \sigma }{B} \big ] \\
			\dot{B} = \epsilon g_3 & f_2 = -\xi & g_3 = - \big [ B \sin^2(\phi) + \beta \cos (\phi) ( B \cos (\phi) + \sigma ) \big ] \\
			\dot{\chi} = \epsilon g_4 & & g_4 = - ( \xi + \sigma \chi )
	\end{array} \end{equation*}
	Hence, $B$ and $\chi$ are constant at $\eps=0$, and for fixed $B$ and $\chi$, the remaining system for $(\xi,\phi)$ possesses the same Hamiltonian structure as \eqref{Lorenz_resc_lim}. Converting to action angle coordinates the system becomes
	\begin{equation*}\begin{split}
			\dot{I} & = \epsilon F_1(I,\theta,B,\chi) \\
			\dot{\theta} & = \Omega(I,B) + \epsilon F_2(I,\theta,B,\chi) \\
			\dot{B} & = \epsilon g_3(I,\theta,B,\chi) \\
			\dot{\chi} & = \epsilon g_4(I,\theta,B,\chi)
	\end{split} \end{equation*}
	where
	\begin{equation*}
		F_1 = \frac{1}{\Omega (I,B)} \big ( \xi g_1 + B \sin (\phi) g_2 - \cos (\phi) g_3 \big ) + \frac{\partial I}{\partial B}|_A  g_3  \hspace{.5 cm} \text{ , } \hspace{.5 cm} 
		F_2 = \frac{\partial \theta}{\partial \xi} g_1 + \frac{\partial \theta}{\partial \phi} g_2 + \frac{\partial \theta}{\partial B} g_3.
	\end{equation*}
	In this case we have a three dimensional Melnikov function given by
	\[ \begin{pmatrix} M_1 \\ M_3 \\ M_4 \end{pmatrix} = \begin{pmatrix} \int_0^t F_1 \big ( I_0,\theta_0 + \Omega(I_0,B_0,\chi_0) s ,B_0,\chi_0 \big ) ds \\ \int_0^t g_3 \big ( I_0,\theta_0 + \Omega(I_0,B_0,\chi_0) s ,B_0,\chi_0 \big ) ds \\ \int_0^t g_4 \big ( I_0,\theta_0 + \Omega(I_0,B_0,\chi_0) s ,B_0,\chi_0 \big )ds \end{pmatrix}. \]
	In order to find the persistent periodic orbits, we need to find the zeros of this vector valued Melnikov function such that the non-degeneracy condition
	\[ \det DM = \det\left(\frac{\partial (M_1,M_3,M_4)}{\partial (I,B, \chi)}\right) \neq  0 \]
	is satisfied.  Since we can write
	\[ M_1 =  \frac{1}{\Omega(I_0,B_0)} \tilde{M}_1 + \frac{\partial I}{\partial B}|_A M_3 \hspace{1 cm} \text{ for } \hspace{1 cm} \tilde{M}_1 = \int_0^T \big [ \xi g_1 + B \sin \phi g_2 - \cos \phi g_3 \big ] dt \]
	it suffices to find $A,B$ such that $\tilde{M}_1 = M_3 = M_4 = 0$.  Explicitly the Melnikov integrals are given via
	\begin{equation*}\begin{split}
			\tilde{M}_1 & = \int_0^T \big [ -\xi^2 \sigma + s \xi \chi + B \beta \cos (\phi) + \beta \sigma \big ] dt  \\
			M_3 & = -\int_0^T \big [ B + (\beta-1) B \cos^2 \phi + \sigma \beta \cos \phi \big ] dt  \\
			M_4 & = \int_0^T ( \xi +  \sigma \chi ) dt.
	\end{split} \end{equation*}
	Since we aim at illustration, we consider $|A| \leq B$ only, and then compute 
	\begin{equation*}\begin{split}
			\tilde{M}_1 & = (B \beta + \beta \sigma ) \frac{4K(k_1)}{\sqrt{B}} - 16 \sigma \sqrt{B} ( E(k_1) - (1-k_1^2 ) K(k_1)) + 8 \beta \sqrt{B} ( E(k_1) - K(k_1)) \\
			M_3 & = -(B + \sigma ) \beta \frac{4K(k_1)}{\sqrt{B}} - 8(2(\beta-1)\sqrt{B} + \frac{\sigma \beta}{\sqrt{B}}) (E(k_1)-K(k_1)) - \frac{16 (\beta-1)\sqrt{B}}{3} ( -2(1+k_1^2)E(k_1) + (2+k_1^2) K(k_1)) \\
			M_4 & = \frac{4K(k_1)}{\sqrt{B}} \sigma \chi_0
	\end{split} \end{equation*}
	As noticed a priori for such an extension of the Lorenz system, the first two Melnikov functions are the same as for \eqref{Lorenz_resc},  the third vanishes if and only if $\chi_0 = 0$, and $M_1$, $M_2$ are independent of $\chi_0$. Hence,  
	\[  DM = \begin{pmatrix}
		\frac{\partial M_1}{\partial I_0} & \frac{\partial M_1}{\partial B_0} & 0 \\ \frac{\partial M_3}{\partial I_0} & \frac{\partial M_3}{\partial B_0} & 0 \\ \frac{\partial M_4}{\partial I_0} & \frac{\partial M_4}{\partial B_0} & \frac{\partial M_4}{\partial \chi_0} \\
	\end{pmatrix} \]
	and the determinant is given by
	\begin{equation*}\begin{split} \mathsf{det} DM  & = \frac{\partial M_4 }{\partial \chi_0} \Big [ \frac{\partial M_1}{\partial I_0} \frac{\partial M_3}{\partial B_0} - \frac{\partial M_1}{\partial B_0} \frac{\partial M_3}{\partial I_0} \Big ] 
			= \frac{4 K(k_1) \sigma }{\sqrt{B}} \Big [ \frac{\partial M_1}{\partial I_0} \frac{\partial M_3}{\partial B_0} - \frac{\partial M_1}{\partial B_0} \frac{\partial M_3}{\partial I_0} \Big ], \end{split} \end{equation*}
	which is non-zero as shown in \S\ref{s:Melnikov}. 
	
	\section{Discussion}
	
	In this paper we have revisited the dynamics of the Lorenz equation in the regime of large Rayleigh number $\rho$, which is known to feature periodic attractors rather than the famous chaotic dynamics \cite{Robbins,Sparrow,costa_knobloch_weiss_1981}. Our main motivation was to study properties of transport of attractors in a parameter regime where states that maximize transport are dynamically unstable. For the Lorenz equations  it was proven in \cite{SouzaDoering1} that maximal transport is realized by the non-zero fixed points,  which are unstable for $\rho>\rho^*, \lambda>1$. However, we found that the literature concerning existence and stability theory of periodic states for large $\rho$ was incomplete. We have therefore provided a rigorous treatment, which essentially confirms the predictions of \cite{Sparrow}. Numerical computations for large finite $\rho$ based on continuation methods and direct simulations have further corroborated these findings. In addition, we have quantified the transport of the periodic attractors and thus the gap of transport compared with the maximum possible. In particular, the transport of the periodic attractors can be arbitrarily small in a parameter range of bistability, where the  states that maximize transport are also stable. 
	Indeed, for fixed $\rho$ we have identified a hysteresis loop in terms of the parameter $\lambda = \frac{\sigma + 1}{\beta + 2}$, which illustrates difficulty to recover from a loss in transport once $\lambda$ exceeds the `tipping point' $\lambda=1$. Moreover, we have computed the stability boundary of periodic attractors in the $(\lambda,\rho)$-plane and found that it extends to relatively low values of $\rho$ below $200$. For fixed $\rho$ we also found a relation to well-known period-doubling bifurcations and symmetric heteroclinic cycles, which produce further regions of bi- and multi-stability of local attractors. 
	
	The Lorenz equations are the crudest mode truncation of the physical model, and there are numerous extensions. For several such generalisations, we have found that our results apply in suitable parameter regimes, in particular for the Lorenz-Stenflo system \cite{LStenflo_1996}. Although our results have no immediate implications in the context of atmospheric convection, we believe they provide a relevant case study for the relation of theoretical bounds and dynamically realized transport.  The approach by perturbing selected solutions from the infinite Rayleigh number limit by exploiting structural properties would be interesting to explore for higher mode truncations and even the viscous Boussinesq equations. Indeed, recent numerical investigations for meaningful bounds in the Boussinesq equation are based on specific solutions and consider  stability properties \cite{wen_goluskin_doering_2022,WenDingChiniKerswell2022}.  We remark that the mode reduced Nusselt number $\mathrm{Nu}=1+\frac{2}{\beta \rho} H(\rho,\beta,\sigma,\textbf{X}_0)$, cf.  \cite{SouzaDoering1}, is bounded by $3$ as $\rho\to\infty$ due to the transport bound from \cite{SouzaDoering1}. However, this is far from the 'ultimate' or `classical' Nusselt number bounds of order $\rho^{1/2}$ or $\rho^{1/3}$ for the PDE model \cite{wen_goluskin_doering_2022,WenDingChiniKerswell2022}.

	The present paper makes a step towards completely settling the question of transport for the Lorenz model.  The set of parameter values for which the transport has not been analytically determined is now reduced to a compact set for which the dynamics are chaotic.  In the large $\rho$ regime, we have analytically determined  stable structures and their transport.  Although we have found numerical evidence for further stable invariant structures, it numerically appears (but remains to be proven) that for fixed $\lambda>1$ and sufficiently large $\rho$ the symmetric periodic orbits are the only attractors.  In the chaotic regime for intermediate Rayleigh numbers the transport is also reduced compared to the non-zero steady states. However, despite the numerous analytical results for the Lorenz attractor, it seems difficult to quantify the transport in that case. It would also be interesting to explore the possible emergence of discrete Lorenz attractors in extended Lorenz systems such as Palmer's \cite{Palmer98}, which is close to a periodic forcing of the Lorenz in a suitable parameter regime.

	\appendix
	
	\section{Positivity of elliptic integral expressions}
	
	\label{app:EllipticIntegralPositivity}
	
	\subsection{First elliptic integral expression}
	
	Here we prove that
	\[ e_2(k) = (4k^2 -1)K + 4(1-2k^2) E > 0 \]
	for any $0 \leq k \leq 1$.  First, note from the explicit formulas
	\[ (4k^2 -1)K + 4(1-2k^2) E = \int_0^1 \frac{1}{\sqrt{(1-t^2)(1-k^2t^2)}} \big [ (4k^2 -1) + 4(1-2k^2) (1-k^2t^2) \big ] dt  \]
	For $0 \leq k^2 \leq 3/4$, the integrand is pointwise positive a.e., hence the integral is positive.  Indeed, for $0 \leq k^2 \leq 1/2$ one has
	\begin{equation*}\begin{split}
			(4k^2 -1) + 4(1-2k^2) (1-k^2t^2) & = 3-4k^2 -4k^2(1-2k^2)t^2 \\ & \geq 3-4k^2 -4k^2(1-2k^2) =  3-8k^2 +8k^4 \\ & > 0
	\end{split} \end{equation*}
	whereas for $1/2 < k^2 \leq 3/4$, one has 
	\[ (4k^2 -1) + 4(1-2k^2) (1-k^2t^2) = 3-4k^2 +4k^2(2k^2-1)t^2 \geq 3-4k^2 \geq 0 \]
	
	On the other hand, for $3/4 < k^2 \leq 1$ the integrand is no longer pointwise positive.  Instead, note that
	\[ \frac{d}{dt} \Big [ \frac{3-4k^2 + 4k^2(2k^2-1)t^2}{\sqrt{(1-t^2)(1-k^2t^2)}} \Big ] = \frac{t \big ( 3 -9k^2+12k^4 - 2k^2(1-2k^2+4k^4)t^2 \big )}{(1-t^2)^{3/2}(1-k^2t^2)^{3/2}} \]
	and since for $k^2 \in [3/4,1]$ it follows that
	\[ 3 -9k^2+12k^4 > 0 \hspace{.5 cm} \text{ , } \hspace{.5 cm} - 2k^2(1-2k^2+4k^4) < 0 \hspace{.5 cm} \text{ , } \hspace{.5 cm} \frac{3 -9k^2+12k^4}{2k^2(1-2k^2+4k^4)} \geq 1 \]
	Hence the integrand is a strictly monotonically increasing function of $t$ for $t \in [0,1]$.  The minimum of the integrand for $t \in [0,1]$ is thus achieved at $t=0$, with minimum equal to $3-4k^2$.  On the other hand note that the integrand tends toward positive infinity as $t \to 1$.  Thus for any $p \geq 0$ we can define $t_p = t_p(k)$ to be the point such that the integrand is equal to $p$, ie
	\begin{equation} \label{EllipticPositivity_TpDef} \frac{3-4k^2 + 4k^2(2k^2-1)t_p^2}{\sqrt{(1-t_p^2)(1-k^2t_p^2)}} = p \end{equation}
	Note that $0 \leq t_{p_1} < t_{p_2} < 1$ for all $0 \leq p_1 < p_2 < \infty$ follows trivially, and furthermore that (\ref{EllipticPositivity_TpDef}) is equivalent to
	\[ a_{k,p} t_p^4 + b_{k,p} t_p^2 + c_{k,p} = 0 \]
	with
	\begin{equation*}\begin{split}
			a_{k,p} & = k^2 ( 16k^2 (2k^2-1)^2 -p^2) \\
			b_{k,p} & = p^2 (1+k^2) - 8k^2 (3-10k^2+8k^4) \\
			c_{k,p} & = (4k^2-3)^2-p^2
	\end{split} \end{equation*}
	In particular, the integrand is less than zero for $t < t_0$ and greater than zero for $t > t_0$, where
	\[ t_0 = \sqrt{\frac{4k^2-3}{4k^2(2k^2-1)} }\]
	which is easily seen to be a monotonically increasing function of $k$ for $k \in [\sqrt{3}/2,1]$.  One can therefore split the domain of integration via
	\[ \Big ( \int_0^{t_0} + \int_{t_0}^1 \Big ) \frac{3-4k^2 + 4k^2(2k^2-1)t^2}{\sqrt{(1-t^2)(1-k^2t^2)}} dt =: I_- + I_+ \]
	and the desired positivity will then follow from lower bounds on $I_-$ and $I_+$.  We consider two cases:
	\begin{enumerate}
		\item $1 \geq k \geq \hat{k}_1$, where $\hat{k}_1 \approx .9266$:  In this case the integrand tends to infinity sufficiently quickly as $t \to 1$ that we can use the most naive bounds on $I_-$.  As mentioned the integrand achieves its minimum at $t =0$, and since $0 < t_0(k) \leq t_0(1) = 1/2$ one has
		\[ I_- = \int_0^{t_0} \frac{3-4k^2 + 4k^2(2k^2-1)t^2}{\sqrt{(1-t^2)(1-k^2t^2)}} dt \geq \int_0^{t_0} \big [ 3-4k^2 \big ] dt \geq -t_0 \geq -1/2 \]
		On the other hand, letting $t_2$ be defined as in (\ref{EllipticPositivity_TpDef}), note
		\[ I_+ \geq \int_{t_2}^1 \frac{3-4k^2 + 4k^2(2k^2-1)t^2}{\sqrt{(1-t^2)(1-k^2t^2)}} dt > \int_{t_2}^1 2dt = 2(1-t_2) \]
		Hence we have $I_+ + I_- > 0$ as long as $t_2 \leq 3/4$.  But, noting that 
		\[ a_{k,2} \big ( \frac{3}{4} \big )^4 + b_{k,2} \big ( \frac{3}{4} \big )^2 + c_{k,2} > 0  \]
		for $k = 1$, it follows from monotonicity of the integrand and the continuity of $t_p(k)$ with respect to $k$ that $t_2 < 3/4$ for all $1 \geq k \geq \hat{k}_1$ where $\hat{k}_1$ is defined such that 
		\[ a_{\hat{k}_1,2} \big ( \frac{3}{4} \big )^4 + b_{\hat{k}_1,2} \big ( \frac{3}{4} \big )^2 + c_{\hat{k}_1,2} = 0  \]
		
		\item $\hat{k}_2 = .927 \geq k \geq \sqrt{\frac{3}{4}}$:  Note this region is deliberately chose to overlap with the region in the first case.  This is easily seen, since for instance one has 
		\[ a_{\hat{k}_2,2} \Big ( \frac{3}{4} \Big )^4 + b_{\hat{k}_2,2} \Big ( \frac{3}{4} \Big )^2 + c_{\hat{k}_2,2} \approx .0026 > 0 \]
		However in this region we have tighter bounds on $I_-$.  Again due the monotonicity the integrand achieves its minimum at $t =0$, and this time $0 < t_0(k) \leq t_0(\hat{k}_2)$, hence one has
		\[ I_- = \int_0^{t_0} \frac{3-4k^2 + 4k^2(2k^2-1)t^2}{\sqrt{(1-t^2)(1-k^2t^2)}} dt \geq \int_0^{t_0} \big [ 3-4\hat{k}_2^2 \big ] dt = (3-4\hat{k}_2^2)t_0(k) \geq (3-4\hat{k}_2^2)t_0(\hat{k}_2) \approx - 0.184 \]
		Again we have
		\[ I_+ > 2(1-t_2) \]
		hence we have $I_+ + I_- > 0$ at least as long as $t_2 \leq 1-\frac{.185}{2} =.9075$.  But it is easily seen that 
		\[ a_{k,2} \big ( .9075 \big )^4 + b_{k,2} \big ( .9075 \big )^2 + c_{k,2} > 0  \]
		for all $.927 \geq k \geq \sqrt{\frac{3}{4}}$, hence the result follows.
	\end{enumerate}
	
	\subsection{Second elliptic integral expression}
	
	\label{app:EllipticIntegralPositivity2}
	
	Here we prove 
	\[ (\frac{73}{20}-2k^2)\frac{K^2}{E^2} - 6 \frac{K}{E} +5 > 0 \]
	for $1 \geq k \geq k_*$, where $k_* \approx .9089$ is defined as the value of $k$ such that $K(k_*) = 2E(k_*)$.  Note that for all $k \leq 1$ one has
	\[ (\frac{73}{20}-2k^2)\frac{K^2}{E^2} - 6 \frac{K}{E} +5 \geq \frac{33}{20}\frac{K^2}{E^2} - 6 \frac{K}{E} +5 \]
	Note that the polynomial
	\[ p_1(x) = \frac{33}{20}x^2 - 6x +5 \]
	is positive for all $x > r_1 \approx 2.34305$.  Note that $K/E$ is monotonically increasing, and $K/E = 2.34305$ when $k = \hat{k}_1 \approx 0.949509$.  On the other hand, when $k \leq 0.94951$, one has
	\[ (\frac{73}{20}-2k^2)\frac{K^2}{E^2} - 6 \frac{K}{E} +5 \geq 1.84686 \frac{K^2}{E^2} - 6 \frac{K}{E} +5 \]
	But the polynomial
	\[ p_2(x) = 1.84686 x^2 - 6x +5 \]
	is non-zero for all $x$, thus proving the bound.

	\subsection{Third elliptic integral expression}
	
	\label{app:EllipticIntegralPositivity3}
	Here we show that
	\begin{equation*}\begin{split} F_2(k) = (2-k^2)E^4 - 8(1-k^2)E^3 K + 6(1-k^2)(2-k^2)E^2K^2 - 2(2-k^2)^2(1-k^2)E K^3 + (2-k^2)(1-k^2)^2 K^4 \end{split}  \end{equation*}
	satisfies $F_2(k)>0$ for all $0 < k < 1$ in two steps, first for $0.998357\approx k_\mathrm{u}<k<1$ and next for $0 < k < k_{\ell} = 0.9984$.
	
	For the interval $k_u < k < 1$, note that the 2nd and 3rd terms of $F_2$ can be rewritten as
	\[ - 8(1-k^2)E^3 K + 6(1-k^2)(2-k^2)E^2K^2 = 2(1-k^2)E^2K(3(2-k^2)K - 4E). \]
	This is strictly positive for all $0 < k < 1$: First,  for $0 < k < \sqrt{\frac{2}{3}}$ one has
	\[ 3(2-k^2)K - 4E > 4(K-E) > 0 \]
	Second, using \textsc{Mathematica}, $\frac{K(\sqrt{2/3})}{E(\sqrt{2/3})} \approx 1.6088<8/3$, and since $K/E$ is monotonically increasing one has 
	\[  3(2-k^2)K - 4E > 3(2-k^2)K - \frac{8}{3}K = (3 (2-k^2) - \frac{8}{3})K > 0. \] 
	Therefore, 
	\begin{equation*}\begin{split} F_2(k) > (2-k^2)E^4 - 2(2-k^2)^2(1-k^2)E K^3 = (2-k^2)E \big ( E^3 - 2(2-k^2)(1-k^2)K^3 \big ). \end{split} \end{equation*}
	Since $(1-k^2)K^3$ equals $\frac{\pi^3}{8}$ at $k=0$, monotonically decreases for $0 < k <1$, and tends to zero as $k\to 1$, there is $k_u$ with
	\[ (1-k_u^2)K(k_u)^3 = \frac{1}{4}. \]
	Using \textsc{Mathematica}, we find $k_u \approx 0.998357$ so that, for $1 > k > k_u$, 
	\[ E^3 - 2(2-k^2)(1-k^2)K^3 \geq E^3 - (1-\frac{k^2}{2}) > 0, \]
	which means $F_2>0$ on this interval.
	
	On the other hand, for $k$ in any neighborhood bounded away from $1$, we can use uniformly convergent series expansions to prove $F_2$ is positive.  First, let 
	\[ p_1 = \Big ( 4 \frac{1-k^2}{2-k^2} \Big )^{1/4} \hspace{1 cm} \text{ , } \hspace{1 cm}  p_1 = \Big ( \frac{(1-k^2)(2-k^2)^3}{4} \Big )^{1/4} \]
	be the positive fourth roots.  One has the identity
	\[ F_2 - (p_1 E - p_2 K)^4 = \frac{k^4}{2-k^2} ( E^4 - (1-k^2)(1-\frac{k^2}{2})^2 K^4 ). \]
	If one defines
	\[ r_+ = \Big ( (1-k^2)(1-\frac{k^2}{2})^2 \Big )^{1/4} \]
	to be the positive fourth root, then one has
	\[ E^4 - (1-k^2)(1-\frac{k^2}{2})^2 K^4 = (E- r_+ K) ( E^3 + r_+ E^2K + r_+^2 EK^2 + r_+^3 K^3 ). \]
	Hence, positivity of $F_2$ is reduced to showing the positivity of the much simpler expression
	\[ E-r_+ K \]
	Since we consider $k \in [0,1-\epsilon]$, we can use the series expansions for $E$ and $K$:
	\begin{equation} K(k) = \frac{\pi}{2} \sum_{n=0}^{\infty} \left ( P_{2n} \right )^2 k^{2n} \hspace{.5 cm} \text{ , } \hspace{.5 cm} E(k) = \frac{\pi}{2} \sum_{n=0}^{\infty} \left ( P_{2n} \right )^2 \frac{k^{2n}}{1-2n}, \qquad 
		P_{2n} = \frac{(2n)!}{2^{2n}(n!)^2}. \label{eq:EK-expansion}\end{equation}
	Since $r_+$ is real analytic on $|k| < 1-\epsilon$ for any $\epsilon > 0$ we can expand
	\[ r_+ = \sum_{n=0}^{\infty} c_n k^{2n}. \]
	One easily finds that $c_0 = 1$, $c_1 = -\frac{1}{2}$, $c_2 = -\frac{1}{16}$ and $c_3 = -\frac{1}{32}$.  Furthermore, the coefficients satisfy the recurrence relation
	\[ (-3+4n)c_n - (8 + 12 n)c_{n+1} + (16+8n)c_{n+2} = 0. \]
	To see this, let $h(x) = \Big ( (1-x)(1-\frac{x}{2})^2 \Big )^{1/4} = \sum_{n=0}^{\infty} c_n x^n$, for which one has
	\[ \frac{h'(x)}{h(x)} =\frac{3x-4}{4x^2-12x+8}. \]
	Hence,
	\[ (4x^2 -12x+8) \Big ( \sum_{n=0}^{\infty} (n+1)c_{n+1} x^n \Big ) = (3x-4) \Big ( \sum_{n=0}^{\infty} c_{n} x^n \Big )   \]
	from which the recursion relation follows.  Next, we claim that $2c_{n+1} \leq c_{n} \leq c_{n+1} < 0$ for all $n \geq 2$.  This follows inductively from the fact that $2c_3\leq c_2 \leq c_3<0$ and that the set $2x \leq y \leq x < 0$ is invariant under the maps that generate the recursion
	\[  \textbf{R}_{n} \begin{pmatrix} x \\ y \end{pmatrix} = \begin{pmatrix} \frac{12n+8}{8n+16} & -\frac{4n-3}{8 n + 16} \\ 1 & 0 \end{pmatrix} \begin{pmatrix} x \\ y \end{pmatrix} \]
	for $n \geq 2$.  To see this, suppose that $2x \leq y \leq x < 0$, $n \geq 2$ and let $(\tilde{x}, \tilde{y} )^T = \textbf{R}_{n} (x,y)^T$.  Then it follows
	\begin{equation*}\begin{split} & \tilde{x} = \frac{12n+8}{8n+16} x - \frac{4n-3}{8n+16}y \leq \frac{n+5}{n+8} x < 0 \\ & \tilde{x} = \frac{12n+8}{8n+16} x - \frac{4n-3}{8n+16}y \geq \frac{8n+11}{8n+16} x \geq x = \tilde{y} \\ & 2\tilde{x} =  \frac{12n+8}{8n+16} 2x - \frac{8n-6}{8n+16}y \leq \frac{8n+28}{8n+16} x \leq x = \tilde{y}.\end{split} \end{equation*}
	Since $c_n<0$ for $n\geq1$, $r_+$ in $(0,1-\epsilon)$ is smaller than any of its finite truncations $r_N:=\sum_{n=0}^Nc_nk^{2n}$.
	Consequently, 
	\begin{equation*}
		E-r_+K\geq E-r_NK =: h_N.
	\end{equation*}
	In particular, $E-r_+K$ is positive whenever $h_N$ is. The functions $h_N$ have an expansion about $k=0$ in $[0,1)$
	\begin{equation*}
		h_N(k)=\frac{\pi}{2}\sum_{n=0}^\infty \tau_n k^{2n}.
	\end{equation*}
	Using the expansion of $E(k)$ and $K(k)$ in \eqref{eq:EK-expansion}, we have
	\begin{equation*}
		\tau_n=\frac{2n}{1-2n}(P_{2n})^2- \sum_{m=1}^{\min(n,N)} c_{m}(P_{2(n-m)})^2.
		\label{eq:tau_n}
	\end{equation*}
	Notice that by choosing any $N\geq 4$, the leading term is always given by $\frac{\pi}{2}\tau_4k^8$, i.e. $\tau_4=\frac{3}{2048}$.  Furthermore, we found numerically that for any choice of $N$ we always had $\tau_n > 0$ for all $n \leq N$.  We therefore obtain a strategy for a proof as follows.  First, considering $N$ large but yet unspecified, split the series into three terms:
	\[ h_N(k)=\frac{\pi}{2} \Big [ \frac{3}{2048} k^8 + \sum_{n=5}^{N-1} \tau_n k^{2n} + \sum_{n=N}^{\infty} \tau_n k^{2n} \Big ] \]
	Since $c_n<0$ for $n\geq1$, we have $\tau_n>\frac{2n}{1-2n}(P_{2n})^2>-(P_{2n})^2$ for $n\geq N$.  One also has monotonicity $P_{2(n+1)}<P_{2n}$, hence
	\begin{equation*}
		\begin{split}
			h_N(k) & \geq \frac{\pi}{2} \Big [ \frac{3}{2048} k^8 + \sum_{n=5}^{N-1} \tau_n k^{2n} - \sum_{n=N}^{\infty} (P_{2n})^2 k^{2n} \Big ] \\ & \geq \frac{\pi}{2} \Big [ \frac{3}{2048} k^8 + \sum_{n=5}^{N-1} \tau_n k^{2n} - k^{2N} \sum_{n=0}^{\infty} (P_{2n})^2 k^{2n} \Big ] \\ & = \frac{3 \pi }{4096} k^8 + \frac{\pi}{2} \sum_{n=5}^{N-1} \tau_n k^{2n} - k^{2N}K(k) .
		\end{split}
	\end{equation*}
	The above holds for all $0 < k < 1$, whereas on the interval $(0,k_{\ell})$ one has
	\begin{equation*}
		h_N(k) \geq k^{8} \Big ( \frac{3 \pi }{4096} - k^{2N-8}K(k_{\ell}) \Big ) + \frac{\pi}{2} \sum_{n=5}^{N-1} \tau_n k^{2n}.
	\end{equation*}
	Note that the first term on the right hand side is always positive for
	\begin{equation}
		\label{IntervalOfPositivity}
		0<k < \min \Big ( \frac{3\pi}{4096\cdot K(k_{\ell})} \Big )^{\frac{1}{2N-8}} .
	\end{equation}
	hence by choosing $N$ sufficiently large we can make this term positive on the entire interval $(0,k_{\ell})$.  One then obtains positivity of $h_N$ on $[0,k_{\ell}]$ if one can verify numerically that for $5 \leq n \leq N-1$ one has $\tau_n > 0$.
	
	By way of example, in the case where $k_{\ell} = \frac{1}{2}$ we find that choosing $N = 9$ gives $\big ( \frac{3\pi}{4096\cdot K(\frac{1}{2})} \big )^{\frac{1}{8}} \approx 0.517$, so this covers the entire interval $(0,\frac{1}{2})$.  Using \textsc{Mathematica}, we find that
	\[ \tau_5 = \frac{21}{8192} \text{ , } \tau_6 = \frac{412}{131072} \text{ , } \tau_7 = \frac{1859}{524288} \text{ , } \tau_8 = \frac{247197}{67108864}. \]
	On the other hand for the interval $k_{\ell} = 0.9984 $ we used \textsc{Mathematica} to find that for $N = 2360$ the interval of positivity in (\ref{IntervalOfPositivity}) includes all of $(0,0.9984)$, and that $\tau_n > 0$ for $5 \leq n \leq N-1$.  
	
	\printbibliography

@article{Agarwal, 
	title={Circuit bounds on stochastic transport in the Lorenz equations}, 
	AUTHOR = {Weady, Scott and Agarwal, Sahil and Wilen, Larry and Wettlaufer, J.S.},
	JOURNAL = {Physics Letters A},
	FJOURNAL = {},
	VOLUME = {382},
	YEAR = {2018},
	NUMBER = {26},
	PAGES = {1731--1737},
	ISSN = {0375-9601},
	URL = {https://www.sciencedirect.com/science/article/pii/S0375960118304493},
	doi = {10.1016/j.physleta.2018.04.035}
	
}

@book{ByrdFriedman,
	title = {Handbook of elliptic integrals for engineers and scientists, 2nd edition.},
	author = {Byrd, P.F. and Friedman, M.D.},
	isbn = {978-3-642-65140-3},
	series = {Grundlehren der mathematischen Wissenschaften},
	year = {1971},
	DOI = {https://doi.org/10.1007/978-3-642-65138-0},
	publisher = {Springer-Verlag New York Heidelberg Berlin,}
}

@article{bykov_1999, 
	title={On systems with separatrix contour containing two saddle-foci}, 
	volume={95}, 
	DOI={},
	journal={J. Math. Sci.}, 
	publisher={}, 
	author={Bykov, V. V.}, 
	year={1999}, 
	pages={2513--2522}}

@article{costa_knobloch_weiss_1981, 
	title={Oscillations in double-diffusive convection}, 
	volume={109}, 
	DOI={10.1017/S0022112081000918},
	journal={Journal of Fluid Mechanics}, 
	publisher={Cambridge University Press}, 
	author={Costa, L. N. Da and Knobloch, E. and Weiss, N. O.}, 
	year={1981}, 
	pages={25--43}}

@article{Curry1978,
	AUTHOR = {Curry, James H.},
	TITLE = {A generalized {L}orenz system},
	JOURNAL = {Comm. Math. Phys.},
	FJOURNAL = {Communications in Mathematical Physics},
	VOLUME = {60},
	YEAR = {1978},
	NUMBER = {3},
	PAGES = {193--204},
	ISSN = {0010-3616},
	MRCLASS = {65P05 (58F15)},
	MRNUMBER = {495037},
	MRREVIEWER = {Georges A. Lebaud},
	URL = {http://projecteuclid.org/euclid.cmp/1103904126},
}

@misc{auto,
	author={E. J. Doedel}, 
	title={AUTO-07P: Continuation and Bifurcation Software for Ordinary Differential Equations},
	url={http://cmvl.cs.concordia.ca/auto}
}

@article{FantuzziGoluskin,
	author = {Fantuzzi, G. and Goluskin, D. and Huang, D. and Chernyshenko, S. I.},
	title = {Bounds for Deterministic and Stochastic Dynamical Systems using Sum-of-Squares Optimization},
	journal = {SIAM Journal on Applied Dynamical Systems},
	volume = {15},
	number = {4},
	pages = {1962-1988},
	year = {2016},
	doi = {10.1137/15M1053347},
	URL = {https://doi.org/10.1137/15M1053347 },
	eprint = {https://doi.org/10.1137/15M1053347}
}

@article{Felicio_2018,
	doi = {10.1088/2399-6528/aaa955},
	url = {https://dx.doi.org/10.1088/2399-6528/aaa955},
	year = {2018},
	month = {02},
	publisher = {IOP Publishing},
	volume = {2},
	number = {2},
	pages = {025028},
	author = {Carolini C Felicio and Paulo C Rech},
	title = {On the dynamics of five- and six-dimensional Lorenz models},
	journal = {Journal of Physics Communications}
}

@article{GlenSpar,
	doi = {10.1007/BF01020649},
	url = {https://doi.org/10.1007/BF01020649},
	year = {1986},
	month = {05},
	volume = {43},
	author = {Glendinning, Paul and Sparrow, Colin},
	title = {T-points: A codimension two heteroclinic bifurcation},
	journal = {Journal of Statistical Physics}
}

@article {Goluskin2018,
	AUTHOR = {Goluskin, David},
	TITLE = {Bounding averages rigorously using semidefinite programming: Mean moments of the {L}orenz system},
	JOURNAL = {J. Nonlinear Sci.},
	FJOURNAL = {Journal of Nonlinear Science},
	VOLUME = {28},
	YEAR = {2018},
	NUMBER = {2},
	PAGES = {621--651},
	ISSN = {0938-8974},
	MRCLASS = {34C29 (37A25 90C22)},
	MRNUMBER = {3770193},
	MRREVIEWER = {M. L. Blank},
	DOI = {10.1007/s00332-017-9421-2},
	URL = {https://doi.org/10.1007/s00332-017-9421-2},
}

@article{Goluskin2021,
	title = {Heat transport bounds for a truncated model of Rayleigh–Bénard convection via polynomial optimization},
	journal = {Physica D: Nonlinear Phenomena},
	volume = {415},
	pages = {132748},
	year = {2021},
	issn = {0167-2789},
	doi = {https://doi.org/10.1016/j.physd.2020.132748},
	url = {https://www.sciencedirect.com/science/article/pii/S0167278920302396},
	author = {Matthew L. Olson and David Goluskin and William W. Schultz and Charles R. Doering}
}

@misc{OlsonDoering2022,
	doi = {10.48550/ARXIV.2203.02067},
	
	url = {https://arxiv.org/abs/2203.02067},
	
	author = {Olson, Matthew L. and Doering, Charles R.},
	
	title = {Heat transport in a hierarchy of reduced-order convection models},
	
	publisher = {arXiv},
	
	year = {2022},
	
	copyright = {arXiv.org perpetual, non-exclusive license}
}

@article{WenDingChiniKerswell2022,
	author = {Wen, Baole  and Ding, Zijing  and Chini, Gregory P.  and Kerswell, Rich R. },
	title = {Heat transport in Rayleigh–Bénard convection with linear marginality},
	journal = {Philosophical Transactions of the Royal Society A: Mathematical, Physical and Engineering Sciences},
	volume = {380},
	number = {2225},
	pages = {20210039},
	year = {2022},
	doi = {10.1098/rsta.2021.0039},
	
	URL = {https://royalsocietypublishing.org/doi/abs/10.1098/rsta.2021.0039},
	eprint = {https://royalsocietypublishing.org/doi/pdf/10.1098/rsta.2021.0039}
}

@article{wen_goluskin_doering_2022, title={Steady Rayleigh–Bénard convection between no-slip boundaries}, volume={933}, DOI={10.1017/jfm.2021.1042}, journal={Journal of Fluid Mechanics}, publisher={Cambridge University Press}, author={Wen, Baole and Goluskin, David and Doering, Charles R.}, year={2022}, pages={R4}}

@article{LiZhang,
	ISSN = {00361399},
	URL = {http://www.jstor.org/stable/2102263},
	author = {Jibin Li and Jianming Zhang},
	journal = {SIAM Journal on Applied Mathematics},
	number = {4},
	pages = {1059--1071},
	publisher = {Society for Industrial and Applied Mathematics},
	title = {New Treatment on Bifurcations of Periodic Solutions and Homoclinic Orbits at High r in the Lorenz Equations},
	urldate = {2022-11-21},
	volume = {53},
	year = {1993}
}

@article{Lorenz,
	title={Deterministic nonperiodic flow},
	author={Lorenz, Edward
	},
	journal={Journal of Atmospheric Sciences},
	volume={20},
	issue={2},
	pages={130-141},
	DOI={https://doi.org/10.1175/1520-0469(1963)020<0130:DNF>2.0.CO;2},
	year={1963},
}

@article{LStenflo_1996,
	doi = {10.1088/0031-8949/53/1/015},
	url = {https://dx.doi.org/10.1088/0031-8949/53/1/015},
	year = {1996},
	month = {01},
	publisher = {},
	volume = {53},
	number = {1},
	pages = {83},
	author = {L Stenflo},
	title = {Generalized Lorenz equations for acoustic-gravity waves in the atmosphere},
	journal = {Physica Scripta}
}

@article{Molteni,
	ISSN = {08948755, 15200442},
	URL = {http://www.jstor.org/stable/26197269},
	author = {Franco Molteni and Laura Ferranti and T. N. Palmer and Pedro Viterbo},
	journal = {Journal of Climate},
	number = {5},
	pages = {777--795},
	publisher = {American Meteorological Society},
	title = {A Dynamical Interpretation of the Global Response to Equatorial Pacific SST Anomalies},
	urldate = {2022-11-11},
	volume = {6},
	year = {1993}
}

@article{Moon2021,
	author = {Moon, Sungju and Baik, Jong-Jin and Hong, Seong-Ho},
	title = {Coexisting Attractors in a Physically Extended Lorenz System},
	journal = {International Journal of Bifurcation and Chaos},
	volume = {31},
	number = {05},
	pages = {2130016},
	year = {2021},
	doi = {10.1142/S0218127421300160},
	URL = {https://doi.org/10.1142/S0218127421300160},
	eprint = {https://doi.org/10.1142/S0218127421300160}
}

@article{Park2021,
	author = {Park,Junho  and Moon,Sungju  and Seo,Jaemyeong Mango  and Baik,Jong-Jin },
	title = {Systematic comparison between the generalized Lorenz equations and DNS in the two-dimensional Rayleigh-Benard convection},
	journal = {Chaos: An Interdisciplinary Journal of Nonlinear Science},
	volume = {31},
	number = {7},
	pages = {073119},
	year = {2021},
	doi = {10.1063/5.0051482},
	URL = { https://doi.org/10.1063/5.0051482},
	eprint = { https://doi.org/10.1063/5.0051482}
}

@article{Palmer98,
	author = "T. N. Palmer",
	title = "Nonlinear Dynamics and Climate Change: Rossby's Legacy",
	journal = "Bulletin of the American Meteorological Society",
	year = "1998",
	publisher = "American Meteorological Society",
	address = "Boston MA, USA",
	volume = "79",
	number = "7",
	doi = "10.1175/1520-0477(1998)079<1411:NDACCR>2.0.CO;2",
	pages=      "1411 - 1423",
	url = "https://journals.ametsoc.org/view/journals/bams/79/7/1520-0477_1998_079_1411_ndaccr_2_0_co_2.xml"
}

@misc{PreprintIvan,
	author = "Ovsyannikov, Ivan",
	title = {Bounds for the Lorenz system with offsets (in preparation)},
	year = {2022}
}

@article{Robbins,
	ISSN = {00361399},
	URL = {http://www.jstor.org/stable/2100965},
	author = {K. A. Robbins},
	journal = {SIAM Journal on Applied Mathematics},
	number = {3},
	pages = {457--472},
	publisher = {Society for Industrial and Applied Mathematics},
	title = {Periodic Solutions and Bifurcation Structure at High R in the Lorenz Model},
	urldate = {2022-11-21},
	volume = {36},
	year = {1979}
}

@article{Shen2015,
	AUTHOR = {Shen, B.-W.},
	TITLE = {Nonlinear feedback in a six-dimensional Lorenz model: impact of an additional heating term},
	JOURNAL = {Nonlinear Processes in Geophysics},
	VOLUME = {22},
	YEAR = {2015},
	NUMBER = {6},
	PAGES = {749--764},
	URL = {https://npg.copernicus.org/articles/22/749/2015/},
	DOI = {10.5194/npg-22-749-2015}
}

@article{Shen2014,
	author = "Bo-Wen Shen",
	title = "Nonlinear Feedback in a Five-Dimensional Lorenz Model",
	journal = "Journal of the Atmospheric Sciences",
	year = "2014",
	publisher = "American Meteorological Society",
	address = "Boston MA, USA",
	volume = "71",
	number = "5",
	doi = "10.1175/JAS-D-13-0223.1",
	pages=      "1701 - 1723",
	url = "https://journals.ametsoc.org/view/journals/atsc/71/5/jas-d-13-0223.1.xml"
}

@article{SouzaDoering1,
	title = {Maximal transport in the Lorenz equations},
	journal = {Physics Letters A},
	volume = {379},
	number = {6},
	pages = {518-523},
	year = {2015},
	issn = {0375-9601},
	doi = {https://doi.org/10.1016/j.physleta.2014.10.050},
	url = {https://www.sciencedirect.com/science/article/pii/S0375960114012067},
	author = {Andre N. Souza and Charles R. Doering}
}

@book{Sparrow,
	title = {The Lorenz Equations: Bifurcations, Chaos, and Strange Attractors},
	author = {Sparrow, Colin},
	isbn = {9780387907758},
	series = {Applied Mathematical Sciences},
	year = {1982},
	DOI = {https://doi.org/10.1007/978-1-4612-5767-7},
	publisher = {Springer}
}

@article{Tucker,
	title = {The Lorenz attractor exists},
	journal = {Comptes Rendus de l'Académie des Sciences - Series I - Mathematics},
	volume = {328},
	number = {12},
	pages = {1197-1202},
	year = {1999},
	issn = {0764-4442},
	doi = {https://doi.org/10.1016/S0764-4442(99)80439-X},
	url = {https://www.sciencedirect.com/science/article/pii/S076444429980439X},
	author = {Warwick Tucker}
}

@article{Wawrzaszek,
	author = {Wawrzaszek, Anna and Krasi\'{n}ska, Agata},
	title = {Hopf Bifurcations, Periodic Windows and Intermittency in the Generalized Lorenz Model},
	journal = {International Journal of Bifurcation and Chaos},
	volume = {29},
	number = {14},
	pages = {1930042},
	year = {2019},
	doi = {10.1142/S0218127419300428},
	URL = { https://doi.org/10.1142/S0218127419300428},
	eprint = { https://doi.org/10.1142/S0218127419300428}
}

@article{WHper,
	author = {Wiggins, Stephen and Holmes, Philip},
	title = {Periodic Orbits in Slowly Varying Oscillators},
	journal = {SIAM Journal on Mathematical Analysis},
	volume = {18},
	number = {3},
	pages = {592-611},
	year = {1987},
	doi = {10.1137/0518046},
	URL = { https://doi.org/10.1137/0518046	},
	eprint = { https://doi.org/10.1137/0518046}
}

@article{WHhom,
	author = {Wiggins, Stephen and Holmes, Philip},
	title = {Homoclinic Orbits in Slowly Varying Oscillators},
	journal = {SIAM Journal on Mathematical Analysis},
	volume = {18},
	number = {3},
	pages = {612-629},
	year = {1987},
	doi = {10.1137/0518047},
	URL = {https://doi.org/10.1137/0518047},
	eprint = { https://doi.org/10.1137/0518047}
}
	
\end{document}